\documentstyle[12pt]{article}
\input{amssym.def}

\newtheorem{thm}{Theorem}
\newtheorem{lem}{Lemma}
\newtheorem{prop}{Proposition}
\newtheorem{rmk}{Remark}

\newtheorem{question}{Question}

\newcounter{exno}
\newenvironment{example}{\begin{quotation}\refstepcounter{exno}{\bf Example \theexno }}{\end{quotation}}
\newenvironment{proof}{{\em Proof: }\par}{$\Box$\medskip}

\newenvironment{caselist}{\left\{\begin{array}{ll}}{\end{array}\right.}

\newcommand{\sri}{strong rank independence }

\newcommand{\deffont}[1]{{\em #1}}

\newcommand{\emphh}[1]{{\em{#1}}}

\newcommand{\refitem}[1]{{\rm(\ref{#1})}}

\newcommand{\eqref}[1]{Equation \ref{eqn:#1}}
\newcommand{\refeq}[1]{\mbox{\rm (\ref{eqn:#1})}}
\newcommand{\secref}[1]{\S\ref{section:#1}}

\newcommand{\quotes}[1]{``#1''}
\newcommand{\resp}[1]{({\em resp.}, #1)}
\newcommand{\ie}{{\em i.e., }}

\newcommand{\eg}{{\em e.g., }}

\newcommand{\alsu}{({almost surely})}

\newcommand{\pr}{^\prime}
\newcommand{\primes}{^{\prime\prime}}
\newcommand{\howmany}{\sharp}
\newcommand{\ilst}[3]{\mbox{$#1=#2,\dots,#3$}}
\newcommand{\mbraces}[1]{\{#1\}}
\newcommand{\braces}[1]{\mbox{$\mbraces{#1}$}}
\newcommand{\bracelst}[2]{\braces{#1,\dots,#2}}
\newcommand{\inblst}[3]{\mbox{$#1\in\bracelst{#2}{#3}$}}
\newcommand{\bracepair}[2]{\braces{#1,#2}}
\newcommand{\pair}[2]{(#1,#2)}
\newcommand{\clint}[2]{\mbox{$[#1,#2]$}}
\newcommand{\unint}{\clint{0}{1}}
\newcommand{\mord}{\prec}
\newcommand{\ord}{\mbox{$\mord$}}
\newcommand{\morder}[2]{#1\mord #2}
\newcommand{\order}[2]{\mbox{$\morder{#1}{#2}$}}
\newcommand{\mReal}{{\Bbb R}}
\newcommand{\Nat}{{\Bbb N}}
\newcommand{\iNat}[1]{\mbox{$#1\in\Nat$}}
\newcommand{\inI}[1]{\mbox{$#1\in I$}}
\newcommand{\inopint}[1]{\mbox{$#1\in (0,1)$}}
\newcommand{\inclint}[1]{\mbox{$#1\in [0,1]$}}
\newcommand{\setsymdiff}{{\scriptstyle{\triangle}}}
\newcommand{\compose}{{\scriptstyle{\circ}}}
\newcommand{\mAeps}{A_{\meps}}

\newcommand{\inAeps}[1]{\mbox{$#1\in\mAeps$}}
\newcommand{\mIeps}{I_{\meps}}
\newcommand{\Ieps}{\mbox{$\mIeps$}}
\newcommand{\inIeps}[1]{\mbox{$#1\in\mIeps$}}

\newcommand{\esssup}{\mbox{ess sup}}

\newcommand{\mcomb}[2]{\left(\begin{array}{c}#1\cr#2\end{array}\right)}
\newcommand{\mloof}[1]{o\left(#1\right)}

\newcommand{\indexlst}[3]{{#1}_{#2},\dots,{#1}_{#3}}
\newcommand{\seq}[3]{\indexlst{#1}{#2}{#3}}
\newcommand{\mtuple}[3]{(\seq{#1}{#2}{#3})}
\newcommand{\tuple}[3]{\mbox{$\mtuple{#1}{#2}{#3}$}}
\newcommand{\setbld}[2]{\{#1 \ :\ #2\}}

\newcommand{\eqdist}{\stackrel{d}{=}}
\newcommand{\eqdef}{:=}
\newcommand{\iid}{{\em i.i.d.}}
\newcommand{\iud}{{\em i.u.d.}}

\newcommand{\map}[3]{\mbox{$#1:#2\to #3$}}
\newcommand{\selfmap}[2]{\map{#1}{#2}{#2}}

\newcommand{\th}[1]{\mbox{${#1}^{th}$}}
\newcommand{\ith}{\th{i}}
\newcommand{\jth}{\th{j}}
\newcommand{\kth}{\th{k}}

\newcommand{\mjs}[1]{j_{#1}}

\newcommand{\id}{\mbox{$id$}}

\newcommand{\msig}{\sigma}
\newcommand{\msigs}[1]{\msig_{#1}}
\newcommand{\msigof}[1]{\msig\left(#1\right)}
\newcommand{\sig}{\mbox{$\msig$}}
\newcommand{\msigbar}{\bar{\msig}}

\newcommand{\ms}{{\frak s}}
\newcommand{\mss}[1]{\ms_{#1}}

\newcommand{\s}{\mbox{$\ms$}}

\newcommand{\mst}{\tilde{\ms}}

\newcommand{\mangset}[2]{{#1}\langle #2 \rangle}

\newcommand{\mstarset}[1]{{#1}^\ast}
\newcommand{\mbXs}{{\mbXto{\ms}}}
\newcommand{\bXs}{\mbox{$\mbXs$}}

\newcommand{\mSym}{{\frak S}}
\newcommand{\Sym}{\mbox{$\mSym$}}
\newcommand{\msiginSym}{{\msig\in\mSym}}
\newcommand{\msinSym}{\ms\in\mSym}
\newcommand{\sinSym}{\mbox{$\msinSym{}$}}
\newcommand{\mspinSym}{\ms\pr\in\mSym}
\newcommand{\spinSym}{\mbox{$\mspinSym{}$}}

\newcommand{\mSyms}[1]{\mSym_{#1}}
\newcommand{\mSymsup}[2]{\mSym{}^{(#1,#2)}}
\newcommand{\Symsup}[2]{\mbox{$\mSymsup{#1}{#2}$}}
\newcommand{\mSymof}[2]{\mSym{} (#1,#2)}

\newcommand{\mstinSym}{\mst\in\mSym}
\newcommand{\stinSym}{\mbox{$\mstinSym{}$}}

\newcommand{\mthet}{\theta}
\newcommand{\thet}{\mbox{$\mthet$}}
\newcommand{\mthets}[1]{\mthet_{#1}}

\newcommand{\mfs}[1]{f_{#1}}
\newcommand{\mfof}[1]{f\left(#1\right)}
\newcommand{\mfsof}[2]{\mfs{#1}\left(#2\right)}

\newcommand{\mfsub}[3]{\mfs{#1,#2,#3}}
\newcommand{\mfsubof}[4]{\mfsub{#1}{#2}{#3}\left(#4\right)}

\newcommand{\mFof}[1]{F\left(#1\right)}

\newcommand{\mfthet}{\mfs{\mthet}}
\newcommand{\fthet}{\mbox{$\mfthet$}}
\newcommand{\mfthetof}[1]{\mfthet\left(#1\right)}

\newcommand{\mfthets}[1]{f_{\mthets{#1}}}

\newcommand{\mfthetsof}[2]{\mfthets{#1}\left(#2\right)}

\newcommand{\mgsub}[3]{g_{#1,#2,#3}}
\newcommand{\mgsubof}[4]{\mgsub{#1}{#2}{#3}\left(#4\right)}

\newcommand{\mdel}{{\delta}}
\newcommand{\mdelof}[1]{\mdel\left(#1\right)}

\newcommand{\meps}{\varepsilon}
\newcommand{\mepsgo}{\meps>0}
\newcommand{\epsgo}{\mbox{$\mepsgo$}}

\newcommand{\mbX}{{\vec{X}}}
\newcommand{\bX}{\mbox{$\mbX$}}
\newcommand{\Xseq}{\mtuple{X}{1}{n}}
\newcommand{\bXseq}{$\mbX=\Xseq$}
\newcommand{\mXs}[1]{X_{#1}}
\newcommand{\Xs}[1]{\mbox{$\mXs{#1}$}}

\newcommand{\mbXto}[1]{{\vec{X}^{#1}}}

\newcommand{\Xtoseq}[1]{\mtuple{X}{{#1}_1}{{#1}_n}}

\newcommand{\mbXsig}{{\mbXto{\msig}}}
\newcommand{\bXsig}{\mbox{$\mbXsig$}}
\newcommand{\Xsigseq}{\Xtoseq{\msig}}

\newcommand{\mbXdown}{{\vec{X}_{\downarrow}}}
\newcommand{\bXdown}{\mbox{$\mbXdown$}}
\newcommand{\Xdownseq}{\mtuple{X}{(1)}{(n)}}
\newcommand{\Xdowntoseq}[1]{\mtuple{X}{({#1}_1)}{({#1}_n)}}

\newcommand{\mXdowns}[1]{X_{(#1)}}
\newcommand{\Xdowns}[1]{\mbox{$\mXdowns{#1}$}}
\newcommand{\mXdone}{\mXdowns{1}}
\newcommand{\Xdone}{\mbox{$\mXdone$}}
\newcommand{\mXdtwo}{\mXdowns{2}}

\newcommand{\mXdthree}{\mXdowns{3}}
\newcommand{\Xdthree}{\mbox{$\mXdthree$}}

\newcommand{\mbXup}{{\vec{X}_{\uparrow}}}
\newcommand{\bXup}{\mbox{$\mbXup$}}

\newcommand{\mbY}{{\vec{Y}}}
\newcommand{\bY}{\mbox{$\mbY$}}
\newcommand{\mYseq}{\mtuple{Y}{1}{n}}
\newcommand{\Yseq}{\tuple{Y}{1}{n}}
\newcommand{\bYseq}{$\mbY=\mYseq$}
\newcommand{\mYs}[1]{Y_{#1}}
\newcommand{\Ys}[1]{\mbox{$\mYs{#1}$}}
\newcommand{\mbYthet}{\mbY_{\mthet}}
\newcommand{\bYthet}{\mbox{$\mbYthet$}}
\newcommand{\mbYdown}{\mbY_{\downarrow}}

\newcommand{\mcalY}{{\cal Y}}
\newcommand{\calY}{\mbox{$\mcalY$}}
\newcommand{\mcalYevent}[1]{\mbraces{\mcalY{}\in #1}}

\newcommand{\mZs}[1]{Z_{#1}}

\newcommand{\mstate}{{\Sigma}}
\newcommand{\state}{\mbox{$\mstate$}}
\newcommand{\mstates}[1]{\mstate_{#1}}
\newcommand{\mstateskl}{\mstates{k,\ell}}

\newcommand{\mprobof}[1]{{\cal P}\left\{#1\right\}}
\newcommand{\probof}[1]{\mbox{$\mprobof{#1}$}}
\newcommand{\mcondprobof}[2]{\mprobof{#1\ |\ #2}}
\newcommand{\condprobof}[2]{\mbox{$\mcondprobof{#1}{#2}$}}

\newcommand{\mrank}{{\frak R}}
\newcommand{\mRs}[1]{\mrank{}_{#1}}
\newcommand{\Rs}[1]{\mbox{$\mRs{#1}$}}
\newcommand{\mrankof}[2]{\mRs{#1}(#2)}
\newcommand{\rankof}[2]{\mbox{$\mrankof{#1}{#2}$}}
\newcommand{\mrelrank}[2]{\mRs{#1,#2}}
\newcommand{\relrank}[2]{\mbox{$\mrelrank{#1}{#2}$}}
\newcommand{\mrelrankof}[3]{\mrelrank{#1}{#2} (#3)}
\newcommand{\relrankof}[3]{\mbox{$\mrelrankof{#1}{#2}{#3}$}}

\newcommand{\mcube}[1]{I^{#1}}
\newcommand{\mncube}{\mcube{n}}
\newcommand{\mkcube}{\mcube{k}}
\newcommand{\ncube}{\mbox{$\mncube$}}
\newcommand{\mdimsimplex}[1]{\mcube{#1}_{\downarrow}}
\newcommand{\dimsimplex}[1]{\mbox{$\mdimsimplex{#1}$}}
\newcommand{\msimplex}{\mdimsimplex{n}}
\newcommand{\simplex}{\mbox{$\msimplex$}}
\newcommand{\mtwosimplex}{\mdimsimplex{2}}
\newcommand{\twosimplex}{\dimsimplex{2}}
\newcommand{\mdiag}{\Delta_n}
\newcommand{\diag}{\mbox{$\mdiag$}}

\newcommand{\mba}{\vec{a}}
\newcommand{\ba}{\mbox{$\mba$}}
\newcommand{\matuple}{\mtuple{a}{1}{n}}

\newcommand{\mas}[1]{a_{#1}}
\newcommand{\as}[1]{\mbox{$\mas{#1}$}}

\newcommand{\projX}[1]{\mbox{\rm proj}_{\downarrow}(#1)}
\newcommand{\projY}[1]{\mbox{\rm proj}_Y(#1)}

\newcommand{\mmof}[1]{m\left(#1\right)}
\newcommand{\mof}[1]{\mbox{$\mmof{#1}$}}

\newcommand{\mLeb}[1]{\mbox{\rm Leb}_{#1}}
\newcommand{\mLebof}[2]{\mLeb{#1}\left(#2\right)}
\newcommand{\mLebnof}[1]{\mLebof{n}{#1}}
\newcommand{\mLeboneof}[1]{\mLebof{1}{#1}}

\newcommand{\mRanks}{{\cal R}}

\newcommand{\mRsup}[2]{\mRanks^{(#1,#2)}}

\newcommand{\mro}{\rho}
\newcommand{\mroof}[1]{\mro{} (#1)}
\newcommand{\roof}[1]{\mbox{$\mroof{#1}$}}
\newcommand{\mros}[2]{\mro{}_{#1,#2}}
\newcommand{\ros}[2]{\mbox{$\mros{#1}{#2}$}}
\newcommand{\mrosof}[3]{\mros{#1}{#2} (#3)}
\newcommand{\rosof}[3]{\mbox{$\mrosof{#1}{#2}{#3}$}}

\newcommand{\kinonen}{\inblst{k}{1}{n}}

\newcommand{\mps}[1]{p_{#1}}

\newcommand{\mPartc}{{\frak P}_c}
\newcommand{\Partc}{\mbox{$\mPartc$}}

\newcommand{\mIs}[1]{I_{#1}}
\newcommand{\Is}[1]{\mbox{$\mIs{#1}$}}

\newcommand{\mlens}[1]{\ell_{#1}}

\newcommand{\mXpart}{{\frak X}}
\newcommand{\mXatom}[2]{\mXpart(#1,#2)}
\newcommand{\sXatom}[2]{${\scriptscriptstyle{\mXatom{#1}{#2}}}$}
\newcommand{\Xatom}[2]{\mbox{$\mXatom{#1}{#2}$}}

\newcommand{\mmu}{{\mu}}
\newcommand{\mubox}{\mbox{$\mu$}}
\newcommand{\mbmu}{\bar{\mmu}}
\newcommand{\mmus}[1]{\mmu_{#1}}
\newcommand{\mus}[1]{\mbox{$\mmus{#1}$}}
\newcommand{\mbmus}[1]{\mbmu_{#1}}
\newcommand{\bmus}[1]{\mbox{$\mbmus{#1}$}}
\newcommand{\mmuof}[1]{\mmu\left(#1\right)}

\newcommand{\mbmusof}[2]{\mbmus{#1}\left(#2\right)}

\newcommand{\mufunct}{{\frak u}}
\newcommand{\mufunctof}[1]{\mufunct{} (#1)}

\newcommand{\mufunctsof}[2]{\mufunct_{#1} (#2)}

\newcommand{\mmuinv}{{\bar{\mu}}}

\newcommand{\mmuinvs}[1]{\mmuinv_{#1}}

\newcommand{\calF}{{\cal F}}

\newcommand{\mLs}[1]{L_{#1}}
\newcommand{\Ls}[1]{\mbox{$\mLs{#1}$}}

\newcommand{\mBs}[1]{B_{#1}}
\newcommand{\Bs}[1]{\mbox{$\mBs{#1}$}}
\newcommand{\mBcs}[1]{B^c_{#1}}
\newcommand{\Bcs}[1]{\mbox{$\mBcs{#1}$}}
\newcommand{\mBt}{\mBs{t}}
\newcommand{\Bt}{\Bs{t}}
\newcommand{\mBtp}{\mBs{t\pr}}

\newcommand{\mBtpp}{\mBs{t\primes}}

\newcommand{\mBct}{\mBcs{t}}
\newcommand{\Bct}{\Bcs{t}}
\newcommand{\mBctp}{\mBcs{t\pr}}
\newcommand{\Bctp}{\Bcs{t\pr}}
\newcommand{\mBctpp}{\mBcs{t\primes}}

\newcommand{\mBfc}[1]{\mBcs{\mfof{#1}}}

\newcommand{\mCs}[1]{C_{#1}}
\newcommand{\Cs}[1]{\mbox{$\mCs{#1}$}}

\newcommand{\mDs}[1]{D_{#1}}
\newcommand{\Ds}[1]{\mbox{$\mDs{#1}$}}

\newcommand{\calN}{{\cal N}}
\newcommand{\mpproc}[2]{\calN[#1,#2]}
\newcommand{\pproc}[2]{\mbox{$\mpproc{#1}{#2}$}}
\newcommand{\mpprocmB}{\mpproc{m}{B}}
\newcommand{\pprocmB}{\mbox{$\mpprocmB$}}
\newcommand{\mpprocmBof}[1]{\mpproc{m}{B}\left(#1\right)}

\newcommand{\mprocof}[1]{\calN\left(#1\right)}
\newcommand{\procof}[1]{\mbox{$\mprocof{#1}$}}
\newcommand{\mcalNp}{{\cal N\pr}}
\newcommand{\calNp}{\mbox{$\mcalNp$}}
\newcommand{\mcalNpp}{{\cal N\primes}}
\newcommand{\calNpp}{\mbox{$\mcalNpp$}}

\newcommand{\mcalNppof}[1]{\mcalNpp{}\left(#1\right)}

\newcommand{\mcalNt}{\calN{}_t}
\newcommand{\calNtof}[1]{\mcalNt{} (#1)}

\newcommand{\moeps}{o (\meps)}

\newcommand{\mgam}{\gamma}
\newcommand{\mgamof}[1]{\mgam\left(#1\right)}
\newcommand{\mgams}[1]{\mgam_{#1}}
\newcommand{\mgamsof}[2]{\mgams{#1}\left(#2\right)}
\newcommand{\eps}{\mbox{$\meps$}}

\newcommand{\mns}[1]{n_{#1}}
\newcommand{\mNs}[1]{N_{#1}}
\newcommand{\mms}[1]{m_{#1}}
\newcommand{\msub}[1]{\mbox{$\mms{#1}$}}
\newcommand{\mMs}[1]{M_{#1}}
\newcommand{\Msub}[1]{\mbox{$\mMs{#1}$}}

\newcommand{\axes}[2]{
                        \put(0,0){\vector(1,0){115}}
                        \put(120,0){\mbox{$#1$}}
                        \put(0,0){\vector(0,1){115}}
                        \put(0,120){\mbox{$#2$}}
}       

\newcommand{\axisunits}{
                        \put(-10,-10){\mbox{$0$}}
                        \put(-10,100){\mbox{$1$}}
                        \put(100,-10){\mbox{$1$}}
        }

\newcommand{\fthetbox}
{                       \thinlines
                                \axes{t}{\mfthetof{t}}
                                \put(0,100){\line(1,0){100}}
                                \put(100,0){\line(0,1){100}}
                        \thicklines
                                \put(60,0){\line(-3,5){60}}
                                \put(60,0){\line(2,5){40}}
        }

\title{Rank Independence\\and\\ Rearrangements of Random Variables
\ }
\author{Alexander Gnedin\\
Utrecht University
\and Zbigniew Nitecki\thanks{Research supported in part by DFG through SFB 170, ``Geometrie und Analysis'' 
at G\"ottingen, and by a travel grant from Tufts University.} \\
Tufts University
\\}

\begin{document}
\maketitle

\begin{abstract}
We study rearrangements $\mtuple{Y}{1}{n}=\mtuple{X}{\msigs{1}}{\msigs{n}}$
(where \sig{} is a random permutation) of an \iid{} sequence of random
variables \tuple{X}{1}{n} uniformly distributed on \clint{0}{1};
in particular we consider rearrangements satisfying the \sri{} condition,
that the rank of \Ys{k} among $\seq{Y}{1}{k}$
is independent of the values of $\seq{Y}{1}{k-1}$. Nontrivial examples
of such rearrangements are the \quotes{travellers' processes} defined by
Gnedin and Krengel.  We show that these are the only examples 
when $n=2$, and when certain restrictive assumptions hold for $n\geq3$;
we also construct a new class of examples of such rearrangements for which the 
restrictive assumptions do not hold.
\end{abstract}

\section{Introduction}\label{section:intro}

A sequence {\bXseq} of numbers can be reordered by means of any permutation
$\ms$ (thought of as the map $i\mapsto\mss{i}$, \ilst{i}{1}{n}) to 
obtain the new sequence which we denote
        $$\mbXto{\ms}\eqdef\Xtoseq{\ms}$$
When the sequence {\bX} consists of random points chosen independently 
according to the uniform distribution on the unit interval $I=\unint$,
for any fixed permutation {\s} the process {\bXs} has the same 
distribution as {\bX}.  The situation changes, however, when {\s} is also 
allowed to vary.

We shall call a sequence of random variables {\bYseq} a \deffont{rearrangement}
of {\bX} if there is a random variable \sig, defined on the same probability
space as {\bX} with values in the symmetric group {\Sym}, such that
{\bY} has the same distribution as {\bXsig}:
        $$\Yseq\eqdist\Xsigseq.$$
Of course, the distribution of {\bXsig} is the same as that of {\bX} when 
{\bX} is {\iid} and {\sig} is 
independent of {\bX}, but in general they can be quite different.  Our
definition does not require that the entries of {\bX}
be uniformly distributed on $I$ or even \iid; it
can be applied to any random process {\bX}.  We will focus for the most 
part on {\bX} {\iid} and uniformly distributed (\iud), 
noting that the transformation 
technique can be used to reduce the case of {\bX} any continuously distributed
{\iid} sequence to the {\iud} case.  However, in Lemma \ref{lem:rearr} 
(\secref{ordstat}) it will be useful to use this idea in a 
non-{\iid} setting.

There are three standard rearrangements of {\bX}:  the sequence itself is
identified with the \deffont{trivial} rearrangement ($\msig=\id$), and we
also have the \deffont{descending} \resp{\deffont{ascending}} rearrangements
{\bXdown} \resp{\bXup} obtained by rearranging according to size.  In keeping
with \cite{gnedinkrengel}, where certain applications to games were 
investigated, we shall focus on the maximal order statistic, hence on
the {\em descending} order, which we number largest-to-smallest with indices 
in parentheses:
        $$\mbXdown = \Xdownseq \mbox{  with  }\mXdowns{i}\geq\mXdowns{i+1},
                \quad{}\ilst{i}{1}{n}.$$

We are interested in this paper in the consequences of certain conditions on 
the rank statistics of a rearrangement.  Given {\bY}, we define the 
\deffont{initial ranks} as
        \begin{equation}\label{eqn:irank}
                \mRs{k}=1+\howmany\setbld{i<k}{\mYs{i}>\mYs{k}},
                        \quad\ilst{k}{1}{n}.
        \end{equation}
Of course, $\mRs{1}=1$ and in general $\mRs{k}\in\bracelst{1}{k}$;  note that
{\Rs{k}} is a {\em relative} rank (it gives only the position of {\Ys{k}} 
relative to the earlier elements in \bY ) and measures positions in 
{\em descending} order:  $\mRs{k}=j$ precisely if {\Ys{k}} is the
{\jth} largest of $\seq{Y}{1}{k}$.  Of course, we can ignore ties,
since they have probability zero.

We will investigate rearrangements {\bY} with the property that
        $$\seq{\mrank}{k+1}{n}\mbox{ are independent of }
                \mtuple{Y}{1}{k}\mbox{ for }\ilst{k}{1}{n-1}$$
which we refer to as \deffont{\sri}.  Note that this is
strictly stronger than independence of the initial ranks. For example,
the trivial rearrangement has independent ranks, with each of the $n!$ possible
rank configurations \tuple{\mrank}{1}{n} equally likely, but the 
distribution of {\Rs{k+1}} conditioned on the values {\tuple{X}{1}{k}}
depends in an essential way on how these points divide the interval.  On
the other hand, the ascending and descending rearrangements {\bXdown} and
{\bXup} induce a deterministic sequence of ranks, as does any rearrangement 
obtained by applying a fixed permutation {\s} to either of these, so 
that \sri holds for these rearrangements in a trivial
way.

A nontrivial family of rearrangements with the \sri 
property are the \quotes{travellers' processes} constructed in 
\cite{gnedinkrengel}.  One can describe these as follows:  imagine the \Xs{i}'s
as giving the locations of various cities;  two travellers leave a specified 
interior point of $I$ (corresponding to {\thet} below) travelling in opposite
directions, toward the two endpoints of $I$, with constant speeds adjusted 
so that they will reach their respective endpoints simultaneously.  The 
reordering of {\tuple{X}{1}{n}}  is then given by the order in which
the various cities are reached by one or the other traveller.  Formally,
these processes can be defined as follows:

\begin{example}\label{example:travellers}(\quotes{{\em Travellers' process}}, \cite{gnedinkrengel})
        Pick the parameter $\mthet\in\unint$ and consider the \quotes{V-shaped}
        function \selfmap{\mfthet}{I} (Figure \ref{fig:fthet}) defined by
        $$\mfthetof{x}=
                \begin{caselist}
                        \frac{\mthet-x}{\mthet}& 0\leq x\leq\mthet\\
                        \frac{x-\mthet}{1-\mthet}&\mthet\leq x \leq1.
                \end{caselist}
        $$

        \begin{figure}\label{fig:fthet}

\begin{center}
        \begin{picture}(140,140)(-10,-10)
                \axisunits
                \fthetbox
                \put(60,-10){\mbox{$\theta$}}
        \end{picture}
\end{center}

  \caption{{\fthet} for the travellers' process}
        \end{figure}
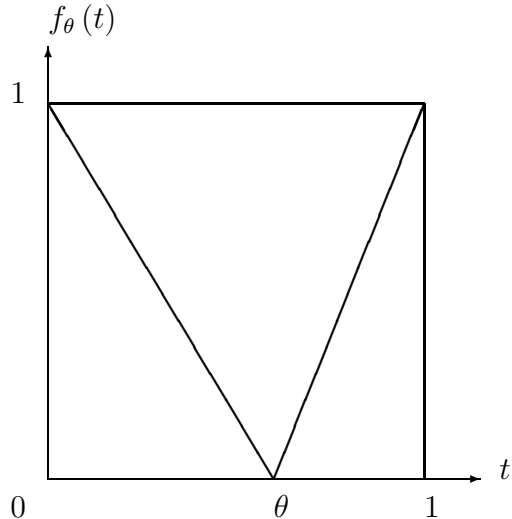

        With probability one, there exists a unique permutation 
        $\msig=\msigof{\mthet,\mbX}\in\Sym$ such that
                $$\mfthetof{\mXs{\msigs{1}}}
                        <\mfthetof{\mXs{\msigs{2}}}
                        <\dots
                        <\mfthetof{\mXs{\msigs{n}}}
                $$
        and the rearrangement
                $$\mbYthet\eqdef\mbXsig$$
        using $\msig=\msigof{\mthet,\mbX}$ has the \sri
        property:  given the values of $\seq{Y}{1}{k}$, we know
        the value of $\mfthetof{\mYs{k}}=\max_{i=1,\dots,k}\mfthetof{\mYs{i}}$
        and that $\mfthetof{\mYs{k+1}}>\mfthetof{\mYs{k}}$;  this tells us
        that {\Ys{k+1}} lies in one of two intervals, the ratio of whose 
        lengths is $\theta/(1-\theta)$, one to the left
        and the other to the right of 
        the interval $\setbld{t}{\mfthetof{t}\leq\mfthetof{\mYs{k}}}$
        ($\mIs{2}$ in Figure \ref{fig:intervals}).  
        One can easily check 
        \cite[section 4]{gnedinkrengel} that in fact the following hold:
        \begin{itemize}
                \item the {\kth} initial rank can only take the extreme values
                        $1$ and $k$;
                      \item the initial rank process $(\mRs{1},\dots,\mRs{n})$
                        can be represented as
                        $$\mRs{k}= J_{k}+ k(1-J_k),\quad\ilst{k}{1}{n}$$
                        where
                        $$J_k={\bf 1}_{[\mthet,1]}
                                \left(\mXs{\msigs{k}}\right)$$
                        are \iid{} Bernoulli variables with 
                        $$\probof{J_k=1}=1-\mthet;$$
                \item $\seq{J}{k+1}{n}$ are independent of 
                        \tuple{Y}{1}{k}.
        \end{itemize}
        This family builds a continuous bridge between the ascending and
        descending rearrangements, with $\mbXup=\mbY_0$ and $\mbXdown=\mbY_1$.
\end{example}

The \sri condition has appeared in various guises in 
connection with different classes of random variables.  While it is
known not to
hold for exchangeable sequences without ties, its relevance to the problems
of Bayesian inference has been discussed in the statistical literature
(see \eg \cite[section 6]{hill1} and\cite{hill2}).  A sequence of independent
(but not identically distributed) random variables satisfying strong rank 
independence was constructed in \cite{hillkennedy}.

In this paper, we investigate the extent to which the \sri
property characterizes the travellers' process of Example 
\ref{example:travellers}.  
\secref{ordstat} gives a framework for thinking about rearrangements
in terms of the descending arrangement {\bXdown}.
In \secref{n2} we will show that when $n=2$, the travellers' 
processes are the only 
rearrangements with the \sri property.  
In \secref{dep}, we show that the \sri{} property forces a certain
dependence between the set of values taken on by the sequence \bX{} and the
rearranging permutation.
In \secref{binary}, we consider the more limited class of {\em binary} 
rearrangements in which the ordering is determined by some real-valued 
attribute (such as the function {\fthet} in Example \ref{example:travellers}) and 
show that for all $n$ the travellers' processes are the only binary 
rearrangements with the \sri property.  In fact, we show
that for binary rearrangements, the strong independence of just a single
rank {\Rs{k}}, $k\in\{2,\dots,n\}$ already forces the rearrangement to be a 
travellers' process.  Finally, in \secref{examples} we discuss some
further examples satisfying the \sri property which
share some features with the travellers' processes and others with the constant
rearrangements in which the components of \bXdown{} are rearranged according
to a fixed element of $\Sym$.

\section{Rerrangements and Order Statistics}\label{section:ordstat}

We shall find it easier to think in terms of the descending rearrangement
{\bXdown} instead of the original sequence {\bX}.  In this section we set
up some machinery to show that this point of view is equivalent to the 
original one.

Note first some general properties of rearrangements.

\begin{lem}\label{lem:rearr}
  Rearrangement is an equivalence relation; that is, for any three
  processes \bX,\bY,$\vec{Z}$ (with the same number of components)
  defined on sufficiently rich probability spaces:
        \begin{enumerate}
        \item\label{rearr:refl}{\bX} is a rearrangement of \bX;
        \item\label{rearr:sym}If {\bY} is a rearrangement of {\bX},
          then {\bX} is a rearrangement of {\bY};
        \item\label{rearr:trans} If {\bY} is a rearrangement of {\bX}
          and $\vec{Z}$ is a rearrangement of {\bY}, then $\vec{Z}$ is
          a rearrangement of {\bX}.
        \end{enumerate}
\end{lem}

\begin{proof}

  \refitem{rearr:refl} is (literally) trivial.

  Note that equality in distribution is preserved by rearrangement, in
  the sense that if
  $$\mbX\eqdist\mbX\pr$$        and $\msiginSym$ is a random
  permutation defined on the same space         as {\bX}, then there
  exists a random permutation $\msig\pr\in\mSym$        defined on the
  same space as $\mbX\pr$ so that
  $$\mbXsig\eqdist\left(\mbX\pr\right)^{\msig\pr}.$$

  Thus, to see \refitem{rearr:sym} we simply note that for any random
  permutation $\msiginSym$, the inverse permutation $\msigbar\in\mSym$
  is also a random permutation, and
  $$\mbX=(\mbXsig)^{\msigbar}.$$

  To see \refitem{rearr:trans}, we note that
  if $\vec{Z}\eqdist\mbY{}^\rho$ for some random $\rho\in\mSym{}$
  (defined on the space for \bY{}) then by \refitem{rearr:sym} there
  is $\rho\pr\in\mSym{}$ (defined on the space for $\vec{Z}$) with
  $\vec{Z}^{\rho\pr}\eqdist\mbY{}$, and hence, since $\mbY{}\eqdist\mbXsig{}$,
  we have $\vec{Z}^{\rho\pr}\eqdist\mbXsig{}$.  But then again we have
  $\rho\primes$ (defined on the space for \bX{}) with $\vec{Z}=
  \left(\vec{Z}^{\rho\pr}\right)^{\bar{\rho}\pr}\eqdist
  \left(\mbXsig{}\right)^{\rho\primes}$.
\end{proof}

We can apply this reasoning in particular to the descending \resp{ascending}
arrangements {\bXdown} \resp{\bXup}.  
Denote by {\simplex} the simplex of descending $n$-tuples in {\ncube}:
        $$\msimplex\eqdef\setbld{\matuple}
                {1\geq a_1\geq a_2\geq\dots \geq a_n\geq0}.$$

There is a \quotes{descending} 
permutation, defined as a map \map{\mdel}{\mncube}{\mSym}, such that for all
$\mba=\matuple\in\mncube$, 
        $$\mba^{\mdelof{\mba}}=\mtuple{a}{\mdel_1}{\mdel_n}\in\msimplex.$$
The value of $\mdel$ is uniquely determined at almost every $\mba\in\mncube$,
specifically, off 
the generalized diagonal \diag{} in {\ncube}:
        $$\mdiag\eqdef\setbld{\matuple\in\mncube}
                {a_i=a_j \mbox{ for some }i\neq j}.$$
Thus, given {\bX} whose entries are continuously distributed on \clint{0}{1}, 
there is a canonical random permutation $\mdel\in\mSym$ defined on the same space as
{\bX} (and uniquely determined a.e.) so that
        $$\mbXdown=\mbXto{\mdel}$$
and hence {\bX} is a rearrangement of 
{\bXdown},
        $$\mbX=\mbXdown^{\bar{\mdel}}.$$
Similar reasoning applies to the ascending rearrangement.
We have, then, as a corollary of Lemma \ref{lem:rearr},

\begin{prop}\label{prop:equiv}
  Given {\bXseq} and {\bYseq} two sequences of random variables as above, the
  following are equivalent:
        \begin{enumerate}
        \item\label{equivi}{\bY} is a rearrangement of {\bX}: for some
          random $\msig\in\mSym$, $\mbY\eqdist\mbXsig;$
        \item\label{equivii}{\bY} and {\bX} have equivalent descending
          rearrangements: $\mbY_\downarrow\eqdist\mbXdown;$
        \item\label{equiviiup}{\bY} and {\bX} have equivalent
          ascending rearrangements: $\mbY_\uparrow\eqdist\mbXup;$
        \item\label{equiviii}{\bY} is a rearrangement of {\bXdown}:
          for some random $\mu\in\mSym$, $\mbY\eqdist\mbXdown^\mu.$
        \end{enumerate}
\end{prop}

The various formulations in Proposition \ref{prop:equiv} can be combined in
a unified picture of rearrangements.
It is easy to see that the \quotes{descending} permutation 
\map{\mdel}{\mncube}{\mSym} is constant on each connected component of 
$\mncube\setminus\mdiag$. Thus, we can identify $\mncube \pmod{0}$ with 
$\msimplex\times\mSym$, by identifying the point 
$\mba\in\mncube\setminus\mdiag$ with the pair 
$\mba^\mdel\in\msimplex,\bar{\mdel}\in\mSym$, where
$\mdel=\mdelof{\mba}$ is the 
\quotes{descending} permutation for \ba, and $\bar{\mdel}$ denotes the 
inverse of $\mdel$ (as a permutation), so that the identification map 
$\msimplex\times\mSym\to\mncube$ is given by $(\mba,\ms)\mapsto\mba^\ms$.

We define a \quotes{state space} 
        $$\mstate\eqdef\msimplex\times\mSym$$
and note that there are two natural \quotes{projections} of {\state}; given 
$(\mba,\ms)\in\mstate$,
        \begin{eqnarray*}
                \projX{\mba,\ms}&\eqdef&\mba\in\msimplex\\
                \projY{\mba,\ms}&\eqdef&\mba^\ms\in\mncube.
        \end{eqnarray*}
Now, if {\bY} is a rearrangement of {\bX}, we can associate to it the 
{\state}-valued random variable
        $$\mcalY\eqdef(\mbXdown,\mu)$$
where $\mu$ is given by Proposition \ref{prop:equiv}\refitem{equiviii}.  We see
that in this case {\bY} and $\mbXdown\eqdist\mbYdown$ can be recovered via the
projections:
        \begin{eqnarray*}
                \mbY&{\eqdist}&\projY{\mcalY}\\
                \mbXdown&{=}&\projX{\mcalY}
        \end{eqnarray*}
Conversely, we have

\begin{lem}\label{lem:state}
        If {\calY} is a random variable with values in 
        $\mstate\eqdef\msimplex\times\mSym$, then 
        $\mbY\eqdef\projY{\mcalY}$ is a
        rearrangement of {\bXdown} (where {\bX} is \iud) if and only if 
        for every measurable set $A\subset\msimplex$,
        \begin{equation}\label{eqn:calY}
                \mprobof{\mcalY\in A\times\mSym}=n!\mLebnof{A}.
        \end{equation}
\end{lem}

\begin{proof}
The right side of the equation is just the normalized Lebesgue measure on 
{\simplex}, or $\probof{\mbXdown{}\in A}$, while the left side is the same
as $\probof{\projX{\mcalY}\in A}$, or equivalently $\probof{\mbYdown\in A}$.
Thus, \eqref{calY} is simply a restatement of the 
requirement that $\mbYdown\eqdist\mbXdown$.
\end{proof}

An advantage of representing a rearrangement \bY{} of \bX{} in terms of 
\bXdown{} and \mubox{} is that it separates data about the values taken
by the variables \Xs{i} from data about their \quotes{arrival times} in
\bY{}.  One can view the \quotes{descending} arrangement \bXdown{} as a
canonical representation of the random (unordered, $n$-point) set of values
\braces{\seq{X}{1}{n}}, and the random permutation \mubox{} as representing the
order in which they are arranged in \bY{}.  \mus{k} gives the \quotes{final}
rank of \Ys{k} among all the variables \Yseq{}, or equivalently \bmus{j} gives
the \quotes{arrival time} for the \jth{} largest value in the sequence
\tuple{Y}{1}{n}.  We shall sometimes refer to \bXdown{} as the \quotes{value
  data} and to \mubox{} as the \quotes{arrival data} for the rearrangement
\bY{}.

An event of the form $ \mtuple{Y}{1}{k}\in{}A\subset{}\mkcube{}$
can be viewed as a
condition on the first $k$ entries of $\projY{\mcalY{}}$;  since
$\projY{\cdot{}}$ is a fixed arrangement of coordinates on each
\quotes{level}
$  \msimplex\times\{\ms\},\quad{}\ms{}\in\Sym{}$ of \state{},
we can formulate the condition as follows:  given $A\subset{}\mcube{k}$
measurable and \sinSym{}, let
\begin{displaymath}
  \mangset{A}{\ms}\eqdef
  \setbld{\mba=\matuple\in\msimplex}
  {\mtuple{a}{{\ms}_1}{{\ms}_k}\in A} \subset\msimplex,
\end{displaymath}
and
 \begin{displaymath}
   \mstarset{A}\eqdef \bigcup_{\msinSym}\mangset{A}{\ms}\times\{\ms\}
   =\mbox{proj}_Y^{-1}\left(A\times\mcube{n-k}\right) \subset\mstate.
 \end{displaymath}
Thus, the event \braces{\mtuple{Y}{1}{k}\in{}A} corresponds in
our representation to \braces{\mcalY{}\in{}\mstarset{A}}.

Rank conditions can also easily be formulated in terms of
$\mcalY{}\in{}\state{}$.  In addition to the initial ranks defined by
\eqref{irank}, we will find it useful to consider other (relative) ranks:  for
any sequence \bYseq{} of random variables without ties, we define $n^2$ \deffont{partial
  ranks} by
\begin{displaymath}
  \mrelrank{j}{k}=\mrelrankof{j}{k}{\mbY{}}\eqdef{}1+\howmany{}
  \setbld{i\leq{}k}{\mYs{i}>\mYs{j}}
  \quad{}\inblst{j,k}{1}{n}.
\end{displaymath}
The initial ranks are given by the special case $j=k$:
\begin{displaymath}
  \mRs{k}=\mrelrank{k}{k},\quad{}\ilst{k}{1}{n};
\end{displaymath}
more generally, for $j\leq{}k$, the numbers \relrank{j}{k} are
\deffont{current ranks}:  if the values of \Yseq{} are displayed consecutively,
then for each \ilst{k}{1}{n} the $k$-tuple $\mtuple{\mrank{}}{1,k}{k,k}$
gives the relative ranking of the first $k$ variables displayed, and encodes
all the rank information known at the \kth{} stage.

To study the interrelationships between the partial ranks more carefully,
we consider their combinatorial analogue, associating to each permutation
\sinSym{} the array \roof{\ms{}} of $n^2$ numbers
\begin{equation}
  \label{eqn:rarray}
  \mrosof{j}{k}{\ms{}}\eqdef{}1+\howmany{}
  \setbld{i\leq{}k}{\mss{i}<\mss{j}}.
\end{equation}
It is clear that, for $\mcalY{}=\pair{\mbXdown{}}{\mmu{}}\in{}\state{}$ as
above the sequence of variables $\mbY{}=\projY{\mcalY{}}$ satisfies
\begin{displaymath}
  \mrelrankof{j}{k}{\mbY{}}=\mrosof{j}{k}{\mu{}}.
\end{displaymath}
The entries on and above the diagonal of \roof{\ms} ($\mrosof{j}{k}{\ms},\ 1
\leq{}j\leq{}k\leq{}n$) will be referred to as the \deffont{upper entries};
they correspond to the current ranks for $\projY{\mcalY{}}$.

\begin{lem}
  \label{lem:rarray} For each \sinSym{}, the numbers $\mros{j}{k}=
  \mrosof{j}{k}{\ms{}}$ defined by \eqref{rarray} satisfy
  \begin{enumerate}
  \item \label{rarrayi} $\mros{j}{k}\in{}\bracelst{1}{k+1}$, the upper entries
    are less than or equal to $k$, and $\mrosof{j}{n}{\ms{}}=\mss{j}$;
  \item \label{rarrayii} the upper values in any column, $\mros{1}{k}\dots{}
    \mros{k}{k}$, are distinct;
  \item \label{rarrayiii} for any $j<k$ with $k>1$,
    \begin{displaymath}
      \mros{k}{k}<\mros{j}{k}\mbox{ iff }\mros{k}{k}\leq{}\mros{j}{k-1};
    \end{displaymath}
  \item \label{rarrayiv} for $j\leq{}k<k\pr$,
    \begin{equation}
      \label{eqn:rarrayiv}
      \mros{j}{k\pr}=\mros{j}{k}+\howmany{}\setbld{\ell{}}{k<\ell\leq{}k\pr
        \mbox{ and }\mros{\ell}{\ell}\leq\mros{j}{\ell-1}}.
    \end{equation}
  \end{enumerate}
\end{lem}

\begin{proof}
  \refitem{rarrayi} is trivial and \refitem{rarrayii} is an immediate
  consequence of the fact that $\mss{1},\dots,\mss{k}$ are distinct.

  To see \refitem{rarrayiii}, note that for any $j,k$ with $k\geq{}2$,
  \begin{equation}
    \label{eqn:nextrank}
    \mros{j}{k}=
    \begin{caselist}
      \mros{j}{k-1}&\mbox{ if }\mss{k}>\mss{j}\cr
      \mros{j}{k-1}+1&\mbox{ if }\mss{k}<\mss{j}
    \end{caselist}
  \end{equation}
  and in either case the two inequalities of \refitem{rarrayiii} are
  equivalent.

  Finally, to see \refitem{rarrayiv}, note that \eqref{rarrayiv} with the
  inequality $\mros{\ell}{\ell}\leq\mros{j}{\ell-1}$ replaced by $\mss{\ell}<
  \mss{j}$ is an easy consequence of \refitem{rarrayii} and the definitions. 
\end{proof}

We can apply \eqref{rarrayiv} recursively to show that any upper entry
\ros{j}{k\pr} ($1\leq{}j\leq{}k\pr\leq{}n$) of $\mro$ is determined uniquely
by any upper entry to its left in the same row (\ros{j}{k}, $k$ fixed,
$j\leq{}k<k\pr$) together with the diagonal entries $\mros{k}{k},
\mros{k+1}{k+1},\dots,\mros{k\pr}{k\pr}$ between.  Conversely, the observation
that the upper entries in column $k$ give the ranking of $\mss{1},\dots,
\mss{k}$ shows that any upper entry \ros{j}{k} ($j\leq{}k$) is also determined
uniquely by the entries in any single column to its right which lie on or
above the same row (\ros{i}{k\pr}, with \ilst{i}{1}{j} and $j\leq{}k<k\pr$
fixed).  To formalize this, for $k\leq{}k\pr$ set
\begin{displaymath}
  \mSymsup{k}{k\pr}\eqdef{}\setbld{(\mss{1},\dots,\mss{k})}
  {\mss{j}\in\bracelst{1}{k\pr}\mbox{ and }\mss{i}\neq\mss{j}\mbox{ for }i\neq{}j}
\end{displaymath}
(so that \Symsup{k}{k} is the set of permutations of \bracelst{1}{k}) and
for $m\leq{}m\pr$ set
\begin{displaymath}
  \mRsup{m}{m\pr}\eqdef{}\setbld{(\indexlst{r}{m}{m\pr})}{r_k\in
    \bracelst{1}{k}\mbox{ for }\ilst{k}{m}{m\pr}}.
\end{displaymath}
Then we have

\begin{rmk}
  \label{rmk:dep} Given $1\leq{}k<k\pr\leq{}n$, there exist functions
  \begin{eqnarray*}
    \mfsub{j}{k}{k\pr}:\mRsup{k}{k\pr}&\to&\bracelst{1}{k\pr}\cr
    \mgsub{j}{k}{k\pr}:\mSymsup{j}{k\pr}&\to&\bracelst{1}{k}
  \end{eqnarray*}
  for $1\leq{}j\leq{}k$ such that for each \sinSym{}, the array
  $\mro{}=\mroof{\ms{}}$ defined by \eqref{rarray} satisfies
  \begin{enumerate}
  \item $\mros{j}{k\pr}=\mfsubof{j}{k}{k\pr}{\mros{j}{k},\mros{k+1}{k+1}
      ,\dots,\mros{k\pr}{k\pr}}$;
  \item $\mros{j}{k}=\mgsubof{j}{k}{k\pr}{\mros{1}{k\pr},\mros{2}{k\pr}
      ,\dots,\mros{j}{k\pr}}$.
  \end{enumerate}
\end{rmk}

Using these functions one easily obtains a bijection for each $k<k\pr$
between $\mSymsup{k}{k}\times\mRsup{k+1}{k\pr}$ (the upper entries in column
$k$ followed by the diagonal through column $k\pr$) and \Symsup{k\pr}{k\pr}
(the upper entries in column $k\pr$).  While an explicit formula for these
bijections is not particularly useful, we will make use of the (well-known)
special case $k=1$, $k\pr=n$.
These bijections also allow us to label the levels of \state{} with
appropriate $n$-tuples of partial ranks, instead of permutations.  
In particular, we can label these levels with initial ranks 
\rosof{k}{k}{\ms{}}.  For $1\leq{}\ell\leq{}k\leq{}n$, let
\begin{displaymath}
  \mSymof{k}{\ell}\eqdef\setbld{\msinSym{}}{\mrosof{k}{k}{\ms}=\ell},
\end{displaymath}
and
\begin{displaymath}
  \mstateskl\eqdef\setbld{(\vec{a},\ms)\in\mstate}
  {\mrelrankof{k}{k}{\projY{(\vec{a},\ms)}}=\ell}=\msimplex{}
  \times\mSymof{k}{\ell}.
\end{displaymath}
Using this notation, we can easily formulate the \sri{} condition in terms
of $\mcalY\in\mstate$.

\begin{rmk}\label{rmk:sri}
  A random variable $\mcalY\in\mstate$ satisfying \eqref{calY}
  in Lemma \ref{lem:state} has the \sri property
  (for $\mbY\eqdist\projY{\mcalY}$) if and only if there exist
  constants
  \begin{displaymath}
    \mps{k,\ell}\geq0,\quad 1\leq\ell\leq k\leq n
  \end{displaymath}
  such that for every $A\subset\mcube{k-1}$,
  \begin{equation}
    \label{eqn:sri}
    \probof{\mcalY\in\mstates{k,\ell}\cap\mstarset{A}}
    =\mps{k,\ell}\probof{\mcalY\in\mstarset{A}},
  \end{equation}
  or equivalently,
  \begin{displaymath}
    \sum_{\ms\in\mSyms{k,\ell}}
    \probof{\mcalY\in\mangset{A}{\ms}\times\{\ms\}}
    =\mps{k,\ell}\probof{\mcalY\in\mstarset{A}}.
  \end{displaymath}
\end{rmk}

Note that for $k=1$, this forces $\mps{1,1}=1$ and for each \inblst{k}{1}{n}
        $$\sum_{\ell=1}^k\mps{k,\ell}=1;$$
of course in general, condition \refeq{sri} is the same as
\begin{equation}
  \label{***}
  \condprobof{\mRs{k}=\ell}{\seq{Y}{1}{k-1}}=\mps{k,\ell}.
\end{equation}

Henceforth, we use this picture to view the descending arrangement 
{\bXdown} as our primary 
object (instead of \bX), using Proposition \ref{prop:equiv}\refitem{equiviii} 
to view any rearrangement
{\bY} as $\mbXdown^\mu$ for some random $\mu\in\mSym$, where
\begin{displaymath}
  \Yseq\eqdist\Xdowntoseq{\mu}.
\end{displaymath}

\section{The case $n=2$}\label{section:n2}

In this section we show that the travellers' processes in Example 
\ref{example:travellers} are the only rearrangements of two (\iud) random 
variables with the \sri condition.  In this case, the combinatorics is 
simplified enormously because {\Sym} contains only two elements, the 
identity $id$ and the transposition $\tau$ ($\tau_1=2,\tau_2=1$).  In terms
of ranks,
        $$\mSymof{2}{1}=\{\tau\},\quad\mSymof{2}{2}=\{id\};$$
it will be convenient to modify our notation from the previous section
slightly and
write for each $A\subset I$
        \begin{eqnarray*}
                \mangset{A}{1}&\eqdef&\mangset{A}{id}
                        =\setbld{\mba\in\mtwosimplex}{a_1\in A}\\
                \mangset{A}{2}&\eqdef&\mangset{A}{\tau}
                        =\setbld{\mba\in\mtwosimplex}{a_2\in A}.
        \end{eqnarray*}

Also, since the \sri condition involves only the two constants 
$\mps{2,1},\mps{2,2}\geq0$ which sum to one, we can express them in terms of a
single parameter $\mthet\in\unint$, with 
        $$\mps{2,1}=1-\mthet,\quad\mps{2,2}=\mthet$$
and the \sri condition is then that for every $A\subset I$,
        $$\probof{\mcalY\in\mangset{A}{2}\times\{\tau\}}
                =\mthet\cdot\probof{\mcalY\in\mstarset{A}}.$$

To simplify our manipulations of certain relations arising from this and 
related conditions, we make the following simple algebraic observation.

\begin{rmk}\label{rmk:equiv}
        Given  $0\leq\mthet\leq1$, let 
        $$\alpha\eqdef\frac{\mthet}{1-\mthet}
                \quad\mbox{(if $\mthet\neq1$)}$$
        and
        $$\alpha^{-1}\eqdef\frac{1-\mthet}{\mthet}
        \quad\mbox{(if $\mthet\neq0$)}.$$
        Then for a given value of {\thet} and any $a,b\geq0$,
        the following are equivalent,
        provided they make sense (\ie $\mthet\neq1$ in \refitem{ratio3} and
        $\mthet\neq0$ in \refitem{ratio4}):
        \begin{enumerate}
                \item\label{ratio1}$a=\mthet(a+b)$;
                \item\label{ratio2}$(1-\mthet)a=\mthet b$;
                \item\label{ratio3}$a=\alpha b$;
                \item\label{ratio4}$\alpha^{-1}a=b$.
        \end{enumerate}
        The same holds if equality is replaced by \quotes{$\leq$}
        in \refitem{ratio1}-\refitem{ratio4}.
\end{rmk}

The travellers' process {\bYthet} from Example \ref{example:travellers} 
(for $n=2$) can be characterized in terms of the function {\fthet}:
        $$\mbY\eqdef(\mYs{1},\mYs{2})\eqdist\mbYthet$$
if and only if
        $$\probof{\mfthetof{\mYs{2}}<\mfthetof{\mYs{1}}}=0$$
(\ie almost surely {\Ys{1}} is the one with the lower \fthet-value).  This
is equivalent to
\begin{equation}
  \label{eqn:concl}
  \forall c\in[0,1)\quad
  \probof{\mfthetof{\mYs{2}}\leq c<\mfthetof{\mYs{1}}}=0
\end{equation}
which is what we will prove.

To this end, fix $c\in[0,1)$ and let
        $$\mPartc\eqdef\{\mIs{1},\mIs{2},\mIs{3}\}$$
be the partition of $I$ into intervals where
        $$\mIs{2}\eqdef\setbld{x}{\mfthetof{x}\leq c}$$
and $\mIs{1},\mIs{3}$ are the components of $\setbld{x}{\mfthetof{x}>c}$, with
$0\in\mIs{1}, 1\in\mIs{3}$ (see Figure \ref{fig:intervals}).  

\begin{figure}\label{fig:intervals}

\begin{center}
        \begin{picture}(140,140)(-20,-20)

                \multiput(0,50)(4,0){25}{\line(1,0){2}}

                \put(-20,25){\makebox(0,0){$c$}}
                \put(-20,30){\vector(0,1){18}}
                \put(-20,20){\vector(0,-1){20}}

                \put(-20,75){\makebox(0,0){$1-c$}}
                \put(-20,80){\vector(0,1){20}}
                \put(-20,70){\vector(0,-1){18}}

                \put(30,-20){\makebox(0,0){$\theta$}}
                \put(25,-20){\vector(-1,0){25}}
                \put(35,-20){\vector(1,0){24}}

                \put(80,-20){\makebox(0,0){$1-\theta$}}
                \put(70,-20){\vector(-1,0){9}}
                \put(90,-20){\vector(1,0){10}}

                \put(15,45){\makebox(0,0){$I_1$}}
                \put(10,45){\vector(-1,0){10}}
                \put(20,45){\vector(1,0){10}}

                \put(90,45){\makebox(0,0){$I_3$}}
                \put(85,45){\vector(-1,0){5}}
                \put(95,45){\vector(1,0){5}}

                \put(50,55){\makebox(0,0){$I_2$}}
                \put(45,55){\vector(-1,0){15}}
                \put(55,55){\vector(1,0){15}}

                \fthetbox

        \end{picture}
\end{center}

        \caption{The partition {\Partc}}
\end{figure}

Denote the length of {\Is{i}} by
        $$\mlens{i}\eqdef\mLeboneof{\mIs{i}}.$$
Using similar triangles in Figure \ref{fig:intervals}, one sees easily that
        \begin{itemize}
                \item $\mlens{2}=c$;
                \item $\mlens{1}/\mthet=\mlens{3}/(1-\mthet)=1-c$
        \end{itemize}
so that in particular (using the notation of Remark \ref{rmk:equiv})
        \begin{equation}\label{eqn:lenratio}
                \mlens{1}=\alpha\mlens{3}.
        \end{equation}

The partition {\Partc} of $I$ gives the product partition 
$\mPartc\times\mPartc$ of $\mcube{2}$, which restricts to {\twosimplex}.
The atoms of this restricted partition are
        $$\mXatom{i}{j}\eqdef
                \setbld{\mba\in\mtwosimplex}{a_1\in\mIs{i},a_2\in\mIs{j}}$$
and since $a_1\geq a_2$ in {\twosimplex}, the only atoms of positive 
measure are the six possibilities for
        $$3\geq i\geq j \geq1$$
(see Figure \ref{fig:Xij}).
\begin{figure}\label{fig:Xij}

\begin{center}

        \begin{picture}(140,140)(-10,-10)
                \axes{\mas{1}}{\mas{2}}
                \axisunits

                \thicklines

                \put(0,0){\line(1,1){100}}

                \put(35,0){\line(0,1){35}}
                \put(65,0){\line(0,1){65}}
                \put(100,0){\line(0,1){100}}

                \put(35,35){\line(1,0){65}}
                \put(65,65){\line(1,0){35}}

                \put(17,17){\makebox(0,0)[tl]{\sXatom{1}{1}}}
                \put(50,17){\makebox(0,0)[tl]{\sXatom{2}{1}}}
                \put(82,17){\makebox(0,0)[tl]{\sXatom{3}{1}}}
                \put(50,50){\makebox(0,0)[tl]{\sXatom{2}{2}}}
                \put(82,50){\makebox(0,0)[tl]{\sXatom{3}{2}}}
                \put(82,82){\makebox(0,0)[tl]{\sXatom{3}{3}}}

        \end{picture}
\end{center}

        \caption{The partition $\mPartc\times\mPartc|\twosimplex$}
\end{figure}

We will also find useful the notation
        $$\mXatom{i}{\ast}\eqdef\bigcup_{j=1}^i\mXatom{i}{j},
        \quad\mXatom{\ast}{j}\eqdef\bigcup_{i=j}^3\mXatom{i}{j}.$$
When $i=j$, \Xatom{i}{j} is a triangle with area
        $$\mLebof{2}{\mXatom{i}{i}}=\frac{\mlens{i}^2}{2}\quad i=1,2,3$$
while for $i>j$, \Xatom{i}{j} is a rectangle, with
        $$\mLebof{2}{\mXatom{i}{j}}=\mlens{i}\mlens{j}\quad 3\geq i>j\geq1.$$
Using the notation  $\mmof{A}\eqdef2\cdot\mLebof{2}{A}$ for the normalized
Lebesgue measure on {\twosimplex}, we see from
\eqref{lenratio} that in particular
        $$\frac{1}{2}\mmof{\mXatom{3}{1}}
                =\alpha^{-1}\mmof{\mXatom{1}{1}}
                =\alpha\mmof{\Xatom{3}{3}},$$
or
        \begin{equation}\label{eqn:arearatio}
                \alpha^{-1}\mmof{\mXatom{1}{1}}+\alpha\mmof{\mXatom{3}{3}}
                        =\mmof{\mXatom{3}{1}}.
        \end{equation}

\begin{thm}\label{thm:n2}
        If a rearrangement $\mbY=(\mYs{1},\mYs{2})$ 
        of $\mbX=(\mXs{1},\mXs{2})$ ({\iud})
        satisfies the \sri condition
        $$\condprobof{\mRs{2}=2}{\mYs{1}\in A}=\mthet
                \quad\mbox{for all $A\subset I$ with }\mLebof{1}{A}>0$$
        then {\bY} is equal in distribution to the corresponding 
        travellers' process of Example \ref{example:travellers}:
        $$\mbY\eqdist\mbYthet.$$
\end{thm}

\begin{proof}
Fix $c\in[0,1)$;  we shall show that the two special cases of the 
hypothesis with $A=\mIs{1}$ \resp{$A=\mIs{3}$}
        \begin{equation}\label{eqn:hyp1}
                \condprobof{\mRs{2}=2}{\mYs{1}\in\mIs{1}}=\mthet
        \end{equation}
        \begin{equation}\label{eqn:hyp2}
                \condprobof{\mRs{2}=2}{\mYs{1}\in\mIs{3}}=\mthet
        \end{equation}
imply \refeq{concl}.

The first hypothesis \refeq{hyp1} can be expressed in terms of {\calY}
and the partition {\Partc} as
        $$\probof{\mcalY\in\mangset{\mIs{1}}{2}\times\{\tau\}}
                =\mthet\probof{\mcalY\in\mstarset{\mIs{1}}}.$$
These sets can be expressed in terms of the partition $\{\mXatom{i}{j}\}$
as follows:
        \begin{eqnarray*}
                \mangset{\mIs{1}}{1}&=&\mXatom{1}{\ast}=\mXatom{1}{1}\\
                \mangset{\mIs{1}}{2}&=&\mXatom{\ast}{1}
                        =\mXatom{1}{1}\cup\mXatom{2}{1}\cup\mXatom{3}{1}\\
                \mstarset{\mIs{1}}&=&\mangset{\mIs{1}}{1}\times\{id\}
                        \cup\mangset{\mIs{1}}{2}\times\{\tau\}.
        \end{eqnarray*}
Using this and rewriting \refeq{hyp1} in form \refitem{ratio2} of 
Remark \ref{rmk:equiv} gives
\begin{eqnarray*}
        \lefteqn{(1-\mthet)
                \probof{\mcalY\in
                \left[\mXatom{1}{1}\cup\mXatom{2}{1}\cup\mXatom{3}{1}\right]
                \times\{\tau\}}}\\
                &=&\mthet\probof{\mcalY\in\mXatom{1}{1}\times\{id\}}
\end{eqnarray*}
and, dropping $\mXatom{1}{1}\times\{\tau\}$ from the event on the left, adding
it to the event on the right, and dividing by $1-\mthet$ gives us
        \begin{equation}\label{eqn:hyp1pr}
                \probof{\mcalY\in
                        \left[\mXatom{2}{1}\cup\mXatom{3}{1}\right]
                        \times \{\tau\}}
                \leq\alpha\probof{\mcalY\in\mXatom{1}{1}\times\mSym}.
        \end{equation}

Similarly, \refeq{hyp2} says
        $$\probof{\mcalY\in\mangset{\mIs{3}}{2}\times\{\tau\}}
                =\mthet\probof{\mcalY\in\mstarset{\mIs{3}}}.$$
This time,
        \begin{eqnarray*}
                \mangset{\mIs{3}}{1}&=&\mXatom{3}{\ast}
                        =\mXatom{3}{3}\cup\mXatom{3}{2}\cup\mXatom{3}{1}\\
                \mangset{\mIs{3}}{2}&=&\mXatom{\ast}{3}
                        =\mXatom{3}{3}\\
                \mstarset{\mIs{3}}&=&\mangset{\mIs{3}}{1}\times\{id\}
                        \cup\mangset{\mIs{3}}{2}\times\{\tau\}
        \end{eqnarray*}
and form \refitem{ratio2} of Remark \ref{rmk:equiv} reads
\begin{eqnarray*}
        \lefteqn{(1-\mthet)
                \probof{\mcalY\in\mXatom{3}{3}\times\{\tau\}}}\\
                &=&\mthet
                \probof{\mcalY\in
                \left[\mXatom{3}{3}\cup\mXatom{3}{2}\cup\mXatom{3}{1}\right]
                \times\{id\}}
\end{eqnarray*}
from which, adding $\mXatom{3}{3}\times\{id\}$ to the event on the left,
dropping it from the right and dividing by {\thet}, we get
        \begin{equation}\label{eqn:hyp2pr}
                \alpha^{-1}\probof{\mcalY\in\mXatom{3}{3}\times\mSym}
                \geq
                \probof{\mcalY\in
                        \left[\mXatom{3}{2}\cup\mXatom{3}{1}\right]
                        \times \{id\}}.
        \end{equation}

        Using Lemma \ref{lem:state} on the right side of \refeq{hyp1pr}
        and the left side of
\refeq{hyp2pr}, writing \refeq{hyp2pr} in reverse order, and 
adding the inequalities gives
\begin{eqnarray*}
        \probof{\mcalY\in\mXatom{2}{1}\times\{\tau\}}
        &+&\probof{\mcalY\in\mXatom{3}{2}\times\{id\}}
        +\probof{\mcalY\in\mXatom{3}{1}\times\mSym}\\
        &\leq&
        \alpha\mmof{\mXatom{1}{1}}
        +\alpha^{-1}\mmof{\mXatom{3}{3}}.
\end{eqnarray*}
But by Lemma \ref{lem:state}, the last term on the left is just 
\mof{\mXatom{3}{1}}, so that \refeq{arearatio} forces
        \begin{equation}\label{eqn:concl2}
        \probof{\mcalY\in\mXatom{2}{1}\times\{\tau\}}
        =\probof{\mcalY\in\mXatom{3}{2}\times\{id\}}
        =0
        \end{equation}

Finally, we analyze \refeq{concl}:
        $$\begin{array}{lcccl}
        \mbraces{\mfthetof{\mYs{2}}\leq c}&=&
                \mcalYevent{\mXatom{\ast}{2}\times\{id\}
                &\cup&
                \mXatom{2}{\ast}\times\{\tau\}}\\
        \{c<\mfthetof{\mYs{1}}\}&=&
                \mcalYevent{\left[\mXatom{1}{\ast}
                        \cup\mXatom{3}{\ast}\right]
                        \times\{id\}
                &\cup&
                \left[\mXatom{\ast}{1}
                        \cup\mXatom{\ast}{3}\right]
                        \times\{\tau\}}
        \end{array}$$
so that
\begin{eqnarray*}
        \{\mfthetof{\mYs{2}}\leq c<\mfthetof{\mYs{1}}\}
        &=&
        \mcalYevent{
                \left[          
                \mXatom{\ast}{2}\cap
                \left(
                        \mXatom{1}{\ast}
                        \cup\mXatom{3}{\ast}
                \right)\right]
              \times\{id\}
              \\
                &\quad&\cup
                \left[
                \mXatom{2}{\ast}\cap
                \left(
                        \mXatom{\ast}{1}
                        \cup\mXatom{\ast}{3}
                \right)\right]
              \times\{\tau\}
              }\\
        &=&
                \mcalYevent{\mXatom{3}{2}\times\{id\}
                \cup\mXatom{2}{1}\times\{\tau\}}
\end{eqnarray*}
and hence \refeq{concl2} is precisely the desired conclusion,
\refeq{concl}.
\end{proof}

\section{Dependence of arrival data on value data}
\label{section:dep}

We saw in \secref{ordstat} that a rearrangement \bY{} is, up to equivalence
in distribution, a function of its \quotes{value data} \bXdown{} and its
\quotes{arrival data} \mubox{}.  It is therefore entirely characterized
by the joint distribution of these data.  One can consider the extent to which
arrival data depends on values; at one extreme the arrival permutation
\mubox{} is independent of \bXdown{}, and at the other it is
\deffont{deterministic} in the sense that for some function
\map{\mufunct{}}{\msimplex{}}{\mSym{}} we have $\mbmu{}
=\mufunctof{\mbXdown{}}$.
Note the distinction between the arrival data \mubox{} and the rearranging
permutation \sig{} (where $\mbXdown^\mmu{}=\mbX{}^\sigma$):
in particular, independence of \mubox{} and \bXdown{} is not equivalent
to independence of \sig{} and \bX{}.
The arrival data for the \quotes{trivial} rearrangement ($\mbXdown^\mmu
=\mbX{}$) is independent of \bXdown{}, since to
recover the \iud{} sequence \bX{} from \bXdown{}, \mubox{} must take each
of the $n!$ possible values in \Sym{} independently of \bXdown{}
with equal probability.
The travellers' processes of \secref{intro} as well as the examples we
will construct in \secref{examples} are by definition deterministic.
Of course, mixed cases are conceivable.

The intersection of the independent and deterministic classes is the set of
\deffont{constant} rearrangements in which \mubox{} takes a single value
in \Sym{} \alsu{}.  A useful \quotes{partial} version of constancy
is that of a \deffont{fixed position}.  The \kth{} position in the
rearrangement \bY{} is \deffont{fixed} at \Xdowns{\ell} if $\mYs{k}
=\mXdowns{\ell}$ \alsu{}.  Clearly, the following are equivalent formulations:
\begin{eqnarray*}
  &\mprobof{\mYs{k}=\mXdowns{\ell}}=1;\cr
  &\mprobof{\mmus{k}=\ell}=1;\cr
  &\mprobof{\mmuinvs{\ell}=k}=1;\cr
  &\mprobof{\mrelrankof{k}{n}{\mbY{}}=\ell}=1.
\end{eqnarray*}
A constant rearrangement is one in which each position is fixed.
\begin{lem}
  \label{lem:fixed}
  Suppose for some rearrangement \bY{} the \kth{} position is fixed and the
  initial ranks \rankof{j}{\bY{}}, \ilst{j}{k}{n} are independent.
  Then each of the partial ranks \relrankof{k}{j}{\bY{}}, \ilst{j}{k}{n}
  is fixed:  that is, it takes a single value \alsu{}.
\end{lem}

\begin{proof}
  In terms of the representation $\mbY{}\eqdist{}\projY{\mbXdown{},\mmu{}}$,
  our hypotheses are that the diagonal entries \ros{j}{j}, \ilst{j}{k}{n} of
  $\mro{}\eqdef\mroof{\mmu{}}$ are independent, and that the last entry
  \ros{k}{n} in the \kth{} row is fixed;  we need to show that then every
  upper entry \ros{k}{j}, \ilst{j}{k}{n} in the \kth{} row is fixed.

  To this end, let \msub{j} \resp{\Msub{j}}, \ilst{j}{k}{n} denote the
  minimum \resp{maximum} of the set $\setbld{r}
  {\mprobof{\mrosof{k}{j}{\mmu}=r}>0}$ of essential values for
  \rosof{k}{j}{\mmu}.
  We claim for \ilst{j}{k}{n-1}
  \begin{equation}
    \label{eqn:claim}
    \mMs{j+1}-\mms{j+1}\geq{}\mMs{j}-\mms{j}.
  \end{equation}
  To see this, note that \eqref{rarrayiv} (Lemma \ref{lem:rarray}) says for
  any \sinSym{} and any \ilst{j}{k}{n-1} that the following analogue of 
  \refeq{nextrank} holds:
  \begin{equation}
    \label{eqn:nextranks}
    \mrosof{k}{j+1}{\ms}=
    \begin{caselist}
      \mrosof{k}{j}{\ms}&\mbox{ if }\mrosof{j+1}{j+1}{\ms}
      >\mrosof{k}{j}{\ms}\cr
      \mrosof{k}{j}{\ms}+1&\mbox{ if }\mrosof{j+1}{j+1}{\ms}
      \leq\mrosof{k}{j}{\ms}.
    \end{caselist}
  \end{equation}
  In particular, $\mMs{j+1}-\mMs{j}$ and $\mms{j+1}-\mms{j}$ are both either
  0 or 1, and \refeq{claim} can fail only if {\it (a)}$\mMs{j+1}=\mMs{j}$
  and {\it (b)} $\mms{j+1}=\mms{j}+1$.  If {\it (b)} occurs, it does so via
  some particular permutation \sinSym{} with $\mprobof{\mmu{}=\ms{}}>0$
  for which 
  \begin{displaymath}
    \mrosof{k}{j}{\ms}=\mms{j}\geq{}\mrosof{j+1}{j+1}{\ms}.
  \end{displaymath}
  We will show that in this case {\it (a)} fails.  Let \spinSym{} with
  $\mprobof{\mmu{}=\ms{}\pr}>0$ and
  \begin{displaymath}
    \mrosof{k}{j}{\ms\pr}=\mMs{j};
  \end{displaymath}
  by independence of initial ranks, there exists \stinSym{} with
  $\mprobof{\mmu{}=\mst}>0$ and
  \begin{displaymath}
    \mrosof{i}{i}{\mst}=
    \begin{caselist}
      \mrosof{i}{i}{\ms\pr}&\mbox{ if }\ilst{i}{k}{j}\cr
      \mrosof{i}{i}{\ms}&\mbox{ if }i=j+1.
    \end{caselist}
  \end{displaymath}
  Recursive application of \eqref{nextranks} gives
  \begin{displaymath}
    \mrosof{k}{j}{\mst}=\mrosof{k}{j}{\ms\pr}=\mMs{j},
  \end{displaymath}
  but then
  \begin{displaymath}
    \mrosof{j+1}{j+1}{\mst}
    =\mrosof{j+1}{j+1}{\ms}
    \leq{}\mms{j}
    \leq{}\mMs{j}
    =\mrosof{k}{j}{\mst}
  \end{displaymath}
  and another application of \eqref{nextranks} then gives $\mMs{j+1}=\mMs{j}+1$,
  contradicting {\it (a)}.
\end{proof}

As we noted in \secref{intro}, the \sri{} condition fails for the trivial
arrangement \bX{} (where \bX{} is \iud{} and \mubox{} is equiprobable and
independent of \bXdown{}),
and holds for any constant arrangement;  this generalizes.

\begin{thm}
  \label{thm:indep}
  If a rearrangement of the \iud{} sequence \bX{} satisfies the \sri{}
  condition and has \bXdown{} and \mubox{} independent, then it is a constant
  rearrangement.
\end{thm}

\begin{proof}
  We will show (by induction) that every position is fixed.
  Pick \kinonen{} and assume that for every $\ell$, $1\leq\ell<k$, the
  \th{\ell} position is fixed.
  We will show that the \kth{} position is fixed.
  The argument rests on three observations.

  The first is that the initial rank \Rs{k} is fixed.  For $k=1$, this is
  trivial.  For $k>1$, our inductive hypothesis, together with Lemma
  \ref{lem:fixed} applied to the \th{\ell} position, \ilst{\ell}{1}{k-1},
  implies that \relrank{k}{\ell} is fixed.  But if $\seq{\mrank}{k,1}{k,k-1}$
  are fixed then Lemma \ref{lem:rarray}\refitem{rarrayii} implies that
  $\mRs{k}=\mrelrank{k}{k}$ is also fixed.

  The second observation is that, for each \inopint{y} and
  $r\in\bracelst{1}{k}$,
  Remark \ref{rmk:dep} and the \sri{} condition formulated as \eqref{sri} (
  Remark \ref{rmk:sri}) give us
  \begin{eqnarray*}
    \mcondprobof{\mmus{k}=r}{\mYs{k}>y}
    &=&\sum_{\mfsubof{k}{k}{n}{\seq{r}{k}{n}}=r}
    \mcondprobof{(\seq{\mrank{}}{k}{n})= (\seq{r}{k}{n})}
    {\mYs{k}>y}\\ &=&\sum_{\mfsubof{k}{k}{n}{\seq{r}{k}{n}}=r}
    \prod_{j=k}^n
    \mps{j,r_j}\cdot\mcondprobof{\mRs{k}=r_k}{\mYs{k}>y}.
  \end{eqnarray*}
  In view of the first observation, the last factor above depends only on
  $r_k$.
  In particular, the conditional probability at the beginning of this
  equation is independent of \inopint{y},
  and hence
  \begin{equation}
    \label{eqn:condprob}
    \mcondprobof{\mmus{k}=r}{\mYs{k}>y}=\mprobof{\mmus{k}=r}.
  \end{equation}

  The third observation is that, if \bXdown{} and \mubox{} are independent,
  we have (again for given $y$ and $r$ as above)
  \begin{displaymath}
    \mprobof{\mmus{k}=r,\mYs{k}>y}
    =\mprobof{\mmus{k}=r,\mXdowns{r}>y}
    =\mprobof{\mmus{k}=r}\cdot\mprobof{\mXdowns{r}>y}
  \end{displaymath}
  and the standard binomial distribution (for \bX{} \iud{}) gives that
  \begin{displaymath}
    \mprobof{\mXdowns{r}>y}
    =\mcomb{n}{r}y^{n-r} (1-y)^r + \mloof{(1-y)^r}\mbox{ as }y\to1.
  \end{displaymath}

  It follows that for each \inblst{r}{1}{k} we have
  \begin{equation}
    \label{eqn:partprob}
    \begin{array}{cll}
      \mprobof{\mmus{k}=r,\mYs{k}>y}
      &=&\mprobof{\mmus{k}=r}\mcomb{n}{r}y^{n-r} (1-y)^r \\
      &&\quad+\mloof{(1-y)^r}\quad\mbox{ as }y\to1
    \end{array}
  \end{equation}
  and, letting $b$ be the minimum value of \mus{k} which appears with positive
  probability,
  \begin{equation}
    \label{eqn:fullprob}
    \begin{array}{cll}
      \mprobof{\mYs{k}>y}
      &=&\sum_{r=b}^k\mprobof{\mmus{k}=r,\mYs{k}>y}\\ 
      &=&\mprobof{\mmus{k}=b}\mcomb{n}{b}y^{n-b} (1-y)^b\\
      &&\quad+\mloof{(1-y)^b}\quad\mbox{ as }y\to1.
    \end{array}
  \end{equation}

  Thus, using \eqref{partprob} with $r=b$ and \eqref{fullprob}, we have
  \begin{displaymath}
    \lim_{y\to1}\mcondprobof{\mmus{k}=b}{\mYs{k}>y}
    =\lim_{y\to1}\frac{\mprobof{\mmus{k}=b,\mYs{k}>y}}{\mprobof{\mYs{k}>y}}=1
  \end{displaymath}
  which, in view of \eqref{condprob}, implies
  \begin{displaymath}
    \mprobof{\mmus{k}=b}=1.
  \end{displaymath}
  Hence position $k$ is fixed.  As \inblst{k}{1}{n} was arbitrary, every
  position is fixed, so the rearrangement is constant and the theorem follows.
\end{proof}

\section{Binary rearrangements}\label{section:binary}

In general, for the deterministic rearrangement given by a function
\map{\mufunct{}}{\msimplex{}}{\mSym{}} (as at the beginning of \secref{dep}),
the position $\mbmusof{i}{\mba{}}$ assigned to the
\ith{} coordinate \as{i} of $\mba\in\msimplex{}$ depends not only on the value
of {\as{i}}, but also on all the other coordinates of {\ba}.  In this section,
we consider those deterministic rearrangements for which the relative
positions assigned to two coordinates depend only on the values of these
two coordinates.  We shall call a map
\map{\mufunct}{\msimplex}{\mSym} \deffont{binary} if there is a subset 
$\calF\subset\mcube{2}$ such that for almost all $\mba\in\msimplex$ and
all $i\neq j$, 
\begin{equation}
  \label{eqn:binmap}
  \mufunctsof{i}{\mba}<\mufunctsof{j}{\mba}\mbox{ iff }
  (\mas{i},\mas{j})\in\calF.
\end{equation}
It is clear that this condition forces $\calF$ to (almost) satisfy the 
basic condition for a total ordering, that for (almost) every pair
$(u,v)\in\mcube{2}$, either $(u,v)\in\calF$
or $(v,u)\in\calF$, but not both.  

We expect a total ordering
to also be transitive.  However, this is not forced by \eqref{binmap}
when $n=2$, as can be seen from the example
\begin{displaymath}
  \mufunct\left(u,v\right)=
  \begin{caselist}
    id &\mbox{if } \frac{1}{3}\leq{}u \leq{}v\leq\frac{2}{3}\\ 
    &\mbox{or } 0\leq{}v<\frac{1}{3} <\frac{2}{3}<u\leq 1\\ 
    \tau\ (\tau_1=2,\tau_2=1)&\mbox{otherwise},
  \end{caselist}
\end{displaymath}
where $(0.2,0.5)$ and $(0.5,0.8)$ but not $(0.2,0.8)$ belong to $\calF$.
(This pathology occurs because our formulation makes {\em every} 
\map{\mufunct}{\mtwosimplex}{\mSym} binary.)
However, for $n\geq3$ transitivity
is forced: if {$(u,v)$} and {$(v,w)$} both belong to $\calF$ and
$\mba\in\msimplex$ has (a permutation of) $(u,v,w)$ as
its first three coordinates, then in $\mba^{\mufunctof{\mba}}$, 
{$u$} \alsu{} precedes {$v$} and {$v$} precedes {$w$}, so {$u$}
precedes {$w$}, hence $(u,w)\in\calF$.  Thus we have

\begin{rmk}\label{rmk:ato}
  For $n\geq3$, every binary map \map{\mufunct}{\msimplex}{\mSym} is
  determined by an \deffont{almost total ordering} of $I$, that is, a binary
  relation {\ord} satisfying:
  \begin{enumerate}
  \item\label{atoi}\emphh{completeness:} the set
    $\setbld{(u,v)\in\mcube{2}} {\mbox{{\em neither
          \order{u}{v} nor \order{v}{u}}}}$
    has measure zero in $\mcube{2}$;

  \item\label{atoii}\emphh{antisymmetry:} the set
    $\setbld{(u,v)\in\mcube{2}} {\mbox{{\em both
          \order{u}{v} and \order{v}{u}}}}$
    has measure zero in $\mcube{2}$;

  \item\label{atoiii}{\em transitivity:} for almost every triple
    $(u,v,w)\in\mcube{3}$ with
    \order{u}{v} and \order{v}{w}, we also
    have \order{u}{w}.
\end{enumerate}
\end{rmk}

One natural way of defining an almost total ordering is by means of a 
measurable function \map{f}{I}{\mReal}, setting {\order{u}{v}} if and 
only if $\mfof{u}<\mfof{v}$:  then properties \refitem{atoi},\refitem{atoii}
and \refitem{atoiii} follow if we assume
$f$ is nonsingular, that is, each level set has measure zero.  Conversely,

\begin{lem}\label{lem:binfunct} 
        Every almost total ordering {\ord} on $I$ is generated by some
        nonsingular measurable function \map{f}{I}{\mReal}, and among
        all such functions (for given ordering \ord{}) there is a unique
        one with values in $I$ which preserves Lebesgue measure.
\end{lem}

\begin{proof}
  Given the almost total ordering {\ord}, define the lower sections
  for \inI{u} by
  \begin{displaymath}
    \mLs{u}\eqdef\setbld{\inI{v}}{v\mord{}u}
  \end{displaymath}
  and set
  \begin{displaymath}
    \mfof{u}\eqdef\mLebof{1}{\mLs{u}}.
  \end{displaymath}

  The transitivity of \ord{} implies that for almost all pairs
  $(u,v)\in I$,
  \begin{equation}\label{eqn:Lusub}
    u\mord v\Rightarrow{}\mLs{u}\subseteq{}\mLs{v}\pmod{0}\quad{}
  \end{equation}
  (that is, \Ls{u} is a subset of $\mLs{u}\cup{}N$ for some null set
  $N$).  In particular, almost surely in $I\times{}I$ we have
  \begin{displaymath}
    \mfof{u}<\mfof{v}\Rightarrow{}u\mord{}v\Rightarrow{}\mfof{u}\leq{}\mfof{v};
  \end{displaymath}
  the second implication is \refeq{Lusub} and the first is its
  contrapositive (with $u$ and $v$ reversed).  It remains to show that
  the set
  \begin{displaymath}
    \setbld{(u,v)}{u\mord{}v\mbox{ and }\mfof{u}=\mfof{v}}
  \end{displaymath}
  has measure zero.

  To this end, pick \inclint{t} and define
  \begin{displaymath}
    \mCs{t}\eqdef\setbld{\inI{u}}{\mfof{u}=t}.
  \end{displaymath}
  We will show that \Cs{t} has measure zero.

  Equation \ref{eqn:Lusub} implies that almost surely,
  \begin{displaymath}
    u\mord{}v,\ u,v\in{}\mCs{t}\Rightarrow{}\mLs{u}=\mLs{v}\pmod{0}.
  \end{displaymath}
  But then completeness of \ord{} insures that for almost every pair $
  (u,v)\in\mCs{t}\times{}\mCs{t}$ we have $\mLs{u}=\mLs{v}\pmod{0}$.
  Fix some $u\pr\in{}\mCs{t}$ such that
  \begin{displaymath}
    \mLs{u}=\mLs{u\pr}\pmod{0}
  \end{displaymath}
  for almost every $u\in{}\mCs{t}$, and let
  \begin{displaymath}
    \mDs{t}\eqdef{}\mCs{t}\cap{}\mLs{u\pr}.
  \end{displaymath}
  Then
  \begin{eqnarray*}
    \mCs{t}\times{}\mDs{t} &=&
    \setbld{(u,v)\in\mCs{t}\times\mCs{t}}{v\mord{}u} \pmod{0}\\ 
    \mDs{t}\times{}\mCs{t} &=&
    \setbld{(u,v)\in\mCs{t}\times\mCs{t}}{u\mord{}v} \pmod{0}
  \end{eqnarray*}
  and by completeness,
  \begin{equation}
    \label{eqn:Dt}
    \mCs{t}\times\mDs{t}\cup\mDs{t}\times\mCs{t}=\mCs{t}\times\mCs{t}\pmod{0}.
  \end{equation}
  But then
  \begin{displaymath}
    \mDs{t}\times\mDs{t}=\mCs{t}\times\mDs{t}\cap\mDs{t}\times\mCs{t}
    =\setbld{(u,v)\in\mCs{t}\times\mCs{t}} {u\mord{}v\mbox{ and
        }v\mord{}u}\pmod{0}
  \end{displaymath}
  must, by antisymmetry, have measure zero.  This implies \Ds{t} has
  measure zero, and hence by the Cavalieri principle, each of the
  (product) sets in Equation \ref{eqn:Dt} has measure zero.  It
  follows that \Cs{t} has measure zero, as required.

  We have shown that
  \begin{displaymath}
    \mLebof{1}{\mCs{t}}=0\mbox{ for each $t$}.
  \end{displaymath}
  First, this implies that
  \begin{displaymath}
    \setbld{(u,v)\in I\times I}{u\mord{}v\mbox{ and
        }\mfof{u}=\mfof{v}=t} \subset{}\mCs{t}\times\mCs{t}
  \end{displaymath}
  has measure zero, so that (by Fubini) almost surely in $I\times I$
  \begin{displaymath}
    u\mord{}v\Rightarrow{}\mfof{u}<\mfof{v}\Rightarrow{}u\mord{}v
  \end{displaymath}
  and second, it implies that $f$ is nonsingular.

  Now, consider the distribution function
  $\mFof{t}\eqdef\mLeboneof{\setbld{u}{\mfof{u}\leq{}t}}$.  Note
  that $F$ is continuous, and for almost every \inI{v}
  \begin{displaymath}
    \mFof{\mfof{v}}=\mLeboneof{\setbld{u}{\mfof{u}<\mfof{v}}}
    =\mLeboneof{\setbld{u}{u\mord{}v}}=\mfof{v}
  \end{displaymath}
  so that $\mFof{t}=t$ for all essential values of $f$; but
  nonsingularity of $f$ implies all values are essential, hence
  $\mFof{t}=t$ $\forall{}\inI{t}$.  This means $f$ is
  measure-preserving.

  Finally, suppose \map{h}{I}{\clint{0}{1}} is another
  measure-preserving function such that almost surely \order{v}{u} iff
  $h (v)<h (u)$; then for all \inI{u},
  \begin{displaymath}
    h(u)=\mLeboneof{\setbld{\inI{v}}{h(v)<h(u)}}
    =\mLeboneof{\mLs{u}}=\mfof{u}.
  \end{displaymath}
\end{proof}

Remark \ref{rmk:ato} and Lemma \ref{lem:binfunct} justify the following 
terminology.  A rearrangement {\bY} of {\bX} is a \deffont{binary 
  rearrangement} if $\mbY\eqdist\mbXdown^{\mu}$, where
$\mbmu=\mufunctof{\mbXdown{}}$ and
\map{\mufunct}{\msimplex}{\mSym} is a binary mapping determined by some almost
total ordering on $I$.
This means (in view of Lemma \ref{lem:binfunct}) that the arrival times
are determined from the values of a (measure-preserving) function
\selfmap{f}{I} via
\begin{displaymath}
  \mfof{\mXdowns{\mmus{1}}}<\mfof{\mXdowns{\mmus{2}}}<\dots<
  \mfof{\mXdowns{\mmus{n}}}.
\end{displaymath}
We will say that the rearrangement is 
\deffont{directed} by $f$, and refer to the family of sets 
        $$\mBt\eqdef\setbld{u}{\mfof{u}\leq t}$$
as the \deffont{filtration} of the rearrangement.

It will be useful for what follows to identify a finite random set with a 
random measure composed of unit point masses.  Suppose $B\subset I$ is a set 
of positive measure.  Let \tuple{U}{1}{m} be an {\iud} sample from $B$; we
define the \deffont{uniform $m$-point process relative to $B$}, {\pprocmB},
by setting, for each Borel set $A\subset I$,
        $$\mpprocmBof{A}\eqdef
                \howmany\setbld{\inblst{i}{1}{m}}{U_i\in A}.$$
Then the familiar formula for multinomial probabilities gives \pprocmB:
if \bracelst{A_1}{A_k} is a (disjoint) partition of $I$, then for each
$k$-tuple \iNat{i_1,\dots,i_k} with $i_1+\dots+i_k= m$,
we have (using $\calN=\mpprocmB$)
\begin{equation}\label{eqn:multi}
        \probof{\procof{A_j}=i_j,\ \ilst{j}{1}{k}}
                =\frac{m!}{i_1!\cdots i_k!}\prod_{j=1}^k
                \left[\frac{\mLeboneof{A_j\cap B}}{\mLeboneof{B}}
                  \right]^{i_j}.
\end{equation}
The uniform point processes determine the original \iud{} samples,
in the sense that given a point process $\calN$ satisfying
\refeq{multi}, we can set up \tuple{U}{1}{m} with
$\mprocof{\mbraces{U_1}}=\dots\mprocof{\mbraces{U_m}}=1$,
$U_1>\dots>U_m$ and then $\mtuple{U}{1}{m}\eqdist\mtuple{X}{(1)}{(m)}$,
where $\seq{X}{1}{m}$ are \iud{} in $B$.

The following properties of {\pprocmB} are straightforward consequences
of \refeq{multi}.

\begin{prop}
  \label{prop:mpoint}
  For any \iNat{m} and $B\subset I$, the uniform  $m$-point processes
  satisfy:
  \begin{enumerate}
  \item\label{mpointi} If $B_i$ are sets converging to $B$ in measure,
    then the processes \pproc{m}{B_i} converge in distribution to
    \pprocmB;
  \item\label{mpointii} If \bracepair{A_1}{A_2} is a partition of $I$,
    then the distribution of the restriction
    $$(\calN=\mpprocmB)|_{A_2}$$ conditioned on $\calN|_{A_1}$, coincides
    with that of
    $$\pproc{m-\mprocof{A_1}}{B\cap A_2}.$$
  \item\label{mpointiii} For $B\pr\subset B$, the distribution of
    $\mpprocmB|_{B\pr}$ conditioned on $\mprocof{B\pr}=m$ coincides with that
    of \pproc{m}{B\pr}.
  \end{enumerate}
\end{prop}

The following relates the processes {\pprocmB} to binary rearrangements.
We use $B^c$ to denote the complement of $B\subset I$ in $I$.

\begin{prop}\label{prop:binpt}
  Suppose {\bY} is a binary rearrangement of {\bX} (\iud) directed by
  $f$, with filtration \braces{\mBt, \inI{t}}.  Then for any
  \ilst{k}{1}{n-1}, the distribution of the random set
  \bracelst{\mYs{k+1}}{\mYs{n}} coincides with that
  of the random point process \pproc{n-k}{\mBfc{\mYs{k}}}.
\end{prop}

\begin{proof}
  Let $\calN=\mpproc{n}{I}$ be the process obtained from {\bX}.  By
  Proposition \ref{prop:mpoint}\refitem{mpointii} for any fixed \inI{t}, the
  distribution of $\calN|_{\mBct}$ conditioned on $\calN|_{\mBt}$ coincides
  with that of \pproc{n-\mprocof{\mBt}}{\mBct}.  This observation
  extends in a straightforward way to the stopping time $T=\mfof{\mYs{k}}$
  which is the moment at which the filtration encounters
  a point of the original process $\calN$ for the \th{k} time.
  By definition, the
  random set \bracelst{\mYs{k+1}}{\mYs{n}} is $\calN|_{\mBcs{T}}$. The
  assertion follows.
\end{proof}

We turn now to binary rearrangements; in view of Theorem \ref{thm:n2} we focus on $n\geq3$.

\begin{prop}\label{prop:binind}
        Suppose {\bY} is a binary rearrangement of {\bX} (\iud), $n\geq3$,
        such that some initial rank \Rs{k}, \inblst{k}{3}{n} is independent
        of  the random variable \tuple{Y}{1}{k-1}.  
        Then almost surely, \Rs{k} takes 
        only its extreme values, $1$ and $k$:
                $$\probof{1<\mRs{k}<k}=0.$$
\end{prop}

\begin{proof}
  By assumption, we have constants $\mps{i}\geq0$, \ilst{i}{1}{k} with
  $\sum\mps{i}=1$ and
  \begin{displaymath}
    \condprobof{\mRs{k}=i}{\indexlst{Y}{1}{k-1}}
    =\mps{i}\quad\ilst{i}{1}{k}.
  \end{displaymath}
  We wish to show
  \begin{displaymath}
    \mps{2}+\dots+\mps{k-1}=0.
  \end{displaymath}

  As usual, we let \selfmap{f}{I} be the function directing {\bY} and
  let
  \begin{displaymath}
    \mBt\eqdef\setbld{u}{\mfof{u}\leq t}, \inI{t}
  \end{displaymath}
  be the associated filtration of {\bY}.  Clearly, {\Bt} is continuous
  in the sense that for all \inI{\alpha, \beta}
  \begin{displaymath}
    \mLeboneof{\mBcs{\alpha}\setsymdiff\mBcs{\beta}} \leq|\alpha-\beta|.
  \end{displaymath}
        
  Fix \inopint{t} and {\epsgo}.  Since $\mLeboneof{\mBt}>0$, we can pick
  an interval {\Ieps} of length {\eps} such that
  \begin{displaymath}
    \mAeps\eqdef\mBt\cap\mIeps
  \end{displaymath}
  has positive measure.

  Momentarily letting $\calN$ denote the $n$-point process defined by
  {\bX}, we note that the event $\mprocof{\mAeps}=k-1,
  \mprocof{\mBct}=n-k+1$ has positive probability (by
  \refeq{multi}) and hence so does the event
  $\indexlst{Y}{1}{k-1}\in\mAeps$ (which is implied by the former).

  Note that if $\indexlst{Y}{1}{k-1}\in\mAeps$ and $1<\mRs{k}<k$, then
  since \Ys{k} lies between the minimum and the maximum of the points
  \bracelst{\mYs{1}}{\mYs{k-1}}, it follows that \inIeps{\mYs{k}}.
  Thus, given \inAeps{\indexlst{Y}{1}{k-1}}, the probability that
  $1<\mRs{k}<k$ is bounded above by the probability that at least one
  of the points $\indexlst{Y}{k}{n}$ belongs to \Ieps.  Now, using
  Proposition \ref{prop:binpt}, let
  \begin{displaymath}
    \calN\eqdist\mpproc{n-k+1}{\mBfc{\mYs{k-1}}}
  \end{displaymath}
  be the ($n-k+1$)-point distribution for $\{\indexlst{Y}{k}{n}\}$.
  Then an easy computation  yields
  \begin{eqnarray*}
    \probof{\mprocof{\mIeps}\geq1}={\cal O} (\meps).
  \end{eqnarray*}
  But
  \begin{displaymath}
    \mps{2}+\dots+\mps{n-1} =\condprobof{1<\mRs{k}<k}
    {\inAeps{\indexlst{Y}{1}{k-1}}} \leq
    \probof{\mprocof{\mIeps}\geq1}
  \end{displaymath}
  and so the proposition follows.
\end{proof}

Using proposition \ref{prop:binind} we can prove the main result of this
section.

\begin{thm}\label{thm:bin}
  Suppose {\bY} is a binary rearrangement of {\bX} (\iud), $n\geq3$,
  such that for some \inblst{k}{2}{n} the initial rank \Rs{k} is
  independent of the random variable $\mtuple{Y}{1}{k-1}$.

  Then {\bY} is equal in distribution to some travellers' process:
  $$\mbY\eqdist\mbYthet.$$
\end{thm}

\begin{proof}
  By Proposition \ref{prop:binind}, \Rs{k} takes only its extreme
  values, $1$ and $k$.  Hence for some \inclint{\mthet} our assumption
  is
  $$\condprobof{\mRs{k}=k}{\indexlst{Y}{1}{k-1}}=\mthet,
  \quad\condprobof{\mRs{k}=1}{\indexlst{Y}{1}{k-1}}=1-\mthet.$$

  As before, we assume {\bY} is directed by the (measure-preserving)
  function $f$ with filtration \Bt, \inI{t}, so that $\mLeboneof{\mBt{}}=t$.

  Fix \inopint{t}, and let $x\pr$ \resp{$x\primes$} be the essential
  infimum \resp{essential supremum} of the set $\mBt\subset I$, and
  set
  $$t\pr\eqdef\lim_{x\downarrow x\pr}
  \esssup\setbld{\mfof{u}}{u\in[x\pr,x]\cap\mBt};$$
  this limit exists because
  $\esssup\setbld{\mfof{u}}{u\in[x\pr,x]}$ decreases with $x$, and
  $t\pr\leq t$.

  For \epsgo, define
  $$\mAeps\pr\eqdef
  \setbld{x}{\mfof{x}\in\clint{t\pr-\meps}{t\pr+\meps}}
  \cap\mBt\cap\clint{x\pr}{x\pr+\meps}.$$ It follows from the
  definition of $x\pr$ that $\mAeps\pr$ has positive measure.

  Similarly, set
  $$t\primes\eqdef\lim_{x\uparrow x\primes}
  \esssup\setbld{\mfof{u}}{u\in[x,x\primes]\cap\mBt},$$ and for \epsgo
  $$\mAeps\primes\eqdef \setbld{x}{\mfof{x}\in
    \clint{t\primes-\meps}{t\primes+\meps}}
  \cap\mBt\cap\clint{x\primes-\meps}{x\primes},$$ so that again
  $t\primes\leq t$ and $\mAeps\primes$ has positive measure.  Note
  that, as $\meps\to0$, we have
        \begin{equation}\label{eqn:23}
          \sup_{x\in\mAeps\pr}|t\pr-\mfof{x}|\to0,\quad
          \sup_{x\in\mAeps\primes}|t\primes-\mfof{x}|\to0.
        \end{equation}
        Now consider the uniform ($n-k+1$)-point process
        $$\mcalNp=\mpproc{n-k+1}{\mBctp}$$ and set $Z\pr$ the
        atom of \calNp{} minimizing $f$.  As in the proof of
        Proposition \ref{prop:binind}, the event
        \braces{\indexlst{Y}{1}{k-1}\in\mAeps\pr} has positive probability.
        Clearly, if $\mtuple{Y}{1}{k-1}\in\mAeps{}\pr$ then
        $\mYs{k}>x\pr+\meps{}$ implies $\mRs{k}=1$, and
        $\mYs{k}<x\pr$ implies $\mRs{k}=k$.
        Thus,
        \begin{eqnarray*}
          \mps{k}&=& \condprobof{\mRs{k}=k}
          {\indexlst{Y}{1}{k-1}\in\mAeps\pr}\\ &=&
          \condprobof{\mRs{k}=k, \mYs{k}\in\clint{x\pr}{x\pr+\meps}}
          {\indexlst{Y}{1}{k-1}\in\mAeps\pr}\\ 
          &&+\condprobof{\mRs{k}=k, \mYs{k}\in\clint{0}{x\pr}}
          {\indexlst{Y}{1}{k-1}\in\mAeps\pr}.
        \end{eqnarray*}

        The first term goes to zero as $\meps\to0$, while by
        Proposition \ref{prop:mpoint}\refitem{mpointiii} and \refeq{23}
        the second converges to \probof{Z\pr\in\clint{0}{x\pr}} so
        $$\mps{k}=\probof{Z\pr\in\clint{0}{x\pr}}.$$

        A similar argument involving conditioning on $\mtuple{Y}{1}{k-1}
        \in\mAeps\primes$ gives
        $$\mps{k}=\probof{Z\primes\in\clint{0}{x\primes}}$$ where
        $Z\primes$ is the atom of
        $\mcalNpp\eqdef\mpproc{n-k+1}{\mBctpp}$ which
        minimizes $f$.

        Next, we claim: $\mBt{}=\clint{x\pr}{x\primes}$ and
        $t=\max{}\mbraces{t\pr,t\primes}$.

        Begin with the case $t\pr\geq{}t\primes$, so that $\mBctpp
        \supset\mBctp$ and $\mBctpp\setminus\mBctp
        =\mBtp\setminus\mBtpp$, and consider the process
        \calNpp{}.  Whenever all atoms of \calNpp{} fall into \Bctp{},
        Proposition \ref{prop:mpoint}\refitem{mpointiii} tells us that
        (conditionally) \calNpp{} agrees in distribution with \calNp{};
        thus,
        \begin{displaymath}
          \mcondprobof{Z\primes\in\clint{0}{x\primes}}
            {\mcalNppof{\mBctpp{}\setminus\mBctp{}}=0}
          =\mprobof{Z\pr\in\clint{0}{x\primes}}.
        \end{displaymath}
        On the other hand, if \calNpp{} has some atoms in
        $\mBctpp{}\setminus\mBctp{}$, then, since $\mfof{x}>t\geq{}t\pr$
        off \clint{x\pr}{x\primes}, we must have
        $Z\primes\in\clint{x\pr}{x\primes}$, and 
        \begin{displaymath}
          \mcondprobof{Z\primes\in\clint{0}{x\primes}}
          {\mcalNppof{\mBctpp{}\setminus\mBctp{}}\geq 1}
          =1.
        \end{displaymath}
        Hence
        \begin{displaymath}
          \begin{array}{lcl}
            \mps{k}&=&\mprobof{Z\primes\in\clint{0}{x\primes}}\\
            &=&\mcondprobof{Z\primes\in\clint{0}{x\primes}}
            {\mcalNppof{\mBctpp{}\setminus\mBctp{}}=0}\cdot
            \mprobof{\mcalNppof{\mBctpp{}\setminus\mBctp{}}=0}\\
           &&\quad +\mcondprobof{Z\primes\in\clint{0}{x\primes}}
           {\mcalNppof{\mBctpp{}\setminus\mBctp{}}\geq1}\cdot
           \mprobof{\mcalNppof{\mBctpp{}\setminus\mBctp{}}\geq1}\\
           &=&\mprobof{Z\pr\in\clint{0}{x\primes}}\cdot
           \mprobof{\mcalNppof{\mBctpp{}\setminus\mBctp{}}=0}\\
           &&\quad+\mprobof{\mcalNppof{\mBctpp{}\setminus\mBctp{}}\geq1}\\
            &\geq&\mprobof{Z\pr\in\clint{0}{x\primes}}
            \geq\mprobof{Z\pr\in\clint{0}{x\pr}}\\
            &=&\mps{k}.
          \end{array}
        \end{displaymath}
        In particular, $\mprobof{Z\pr\in\clint{0}{x\primes}}=
        \mprobof{Z\pr\in\clint{0}{x\pr}}$ implies that
        $\mprobof{Z\pr\in\clint{x\pr}{x\primes}}=0$, a situation possible
        iff \calNp{} \alsu{} puts no atoms in \clint{x\pr}{x\primes},
        or equivalently iff $\mLeboneof{\clint{x\pr}{x\primes}\cap\mBctp{}}=0$,
        which in turn means $\clint{x\pr}{x\primes}\subset\mBtp{}$ $\pmod{0}$
        so that $\mfof{x}\leq{}t\pr$ \alsu{} on \clint{x\pr}{x\primes}.
        Again, since $\mfof{x}\geq{}t$ \alsu{} outside \clint{x\pr}{x\primes}
        and $f$ preserves measure, we must have $t=t\pr$ and, since
        $x\pr,x\primes$ are the essential bounds on \Bt{}, it follows
        also that $\mBt{}=\clint{x\pr}{x\primes}$, and in particular
        $x\primes-x\pr=t$ .

        The argument in case $t\pr\leq{}t\primes$ is similar, involving
        two computations of $\mps{1}$.

        Having established the claim, we now consider each of the endpoints
        of \Bt{} as a non-increasing \resp{non-decreasing} function
        $x\pr (t)$ \resp{$x\primes(t)$}, with
        \begin{displaymath}
          x\primes (t)-x\pr (t)=t\quad\mbox{ for all } t.
        \end{displaymath}
        We wish to compute the derivative of $x\pr (t)$.
        Let $\mcalNt{}$ be a uniform $(n-k+1)$-point process on \Bct{}, and
        $\mZs{t}$ be the atom of $\mcalNt{}$ where $f$ is minimized.
        Arguments like those above give
        \begin{displaymath}
          \begin{array}{lcl}
            \mps{k}&=&\mprobof{\mZs{t}\in\clint{0}{x\pr (t)}}\\
            &=&\mcondprobof{\mZs{t}\in\clint{0}{x\pr (t)}}
            {\calNtof{\mBcs{t+\meps{}}}=n-k+1}\cdot
            \mprobof{\calNtof{\mBcs{t+\meps{}}}=n-k+1}\\
            &&+\mcondprobof{\mZs{t}\in\clint{0}{x\pr (t)}}
            {\calNtof{\mBct{}\setminus\mBcs{t+\meps}}=1}
            \cdot\mprobof{\calNtof{\mBct{}\setminus\mBcs{t+\meps}}=1}\\
            &&+\moeps{}\\
            &=&\mprobof{\mZs{t+\meps}\in\clint{0}{x\pr (t+\meps)}}
            \cdot\mprobof{\calNtof{\mBcs{t+\meps}}=n-k+1}\\
            &&+\mcondprobof{\calNtof{\clint{x\pr (t+\meps)}{x\pr (t)}}=1}
            {\calNtof{\mBct{}\setminus\mBcs{t+\meps}}=1}\cdot
            (n-k+1)\meps (1-t)^{-1}\\
            &&+\moeps{}\\
            &=&\mps{k}\cdot\mprobof{\calNtof{\mBcs{t+\meps}}=n-k+1}\\
            &&+\left(\frac{x\pr(t)-x\pr(t+\meps)}{\meps}\right)
            \cdot (n-k+1)\meps (1-t)^{-1}\\
            &&+\moeps.
          \end{array}
        \end{displaymath}
        Rearranging terms and letting $\meps\to0$ we find that the
        derivative is 
        \begin{displaymath}
          \frac{dx\pr(t)}{dt}=-\mps{k}.
        \end{displaymath}

        This implies
        \begin{displaymath}
          x\pr(t)=-\mps{k}\cdot t+\mps{k},\quad x\primes (t)=\mps{1}\cdot t
          + (1-\mps{1})
        \end{displaymath}
        which in turn forces $f$ to equal \fthet{} with $\thet=\mps{k}$
        (up to a null set).
      \end{proof}

\section{Further examples}\label{section:examples}

So far, the only examples of rearrangements with the \sri property have  been
the travellers' processes of Example \ref{example:travellers} and the constant 
rearrangements $\mbXdown^\ms$, where $\msinSym$ is a fixed permutation.
In this section we construct multiparameter families of deterministic
rearrangements with the \sri property which combine features of both the
travellers' processes and constant rearrangements, but are of neither type.
The idea is that if the position of some $\mXdowns{k}$ in {\bY} is 
fixed, then we can use it to partition $I$ into two subintervals
$\mIs{1}\cup\mIs{2}$, with $n-k$ \resp{$k-1$} points uniformly distributed
on \Is{1} \resp{\Is{2}}, whatever value $\mXdowns{k}$ takes;  these two
\quotes{sub}-processes are independent, and we can rearrange each separately.

Keep in mind the following features of our examples so far:
\begin{itemize}
        \item for a constant rearrangement, each initial rank \Rs{k} almost
                surely takes a single value;
        \item for the travellers' process (or by theorem \ref{thm:bin},
                any binary rearrangement),
                each initial rank \Rs{k} takes only the extreme values
                $1$ and $k$.
\end{itemize}

Before giving a general construction, we consider two specific examples:

\begin{example}\label{example:2}
Take $n=3$ and choose \inopint{\mthet}.  Now set
        $$\mYs{1}=\mXdone,$$
and given {\Xdone}, let \map{\gamma}{\clint{0}{\mXdone}}{I} be the 
unique linear, order-preserving bijection, $t\mapsto t/\mXdone$.
(Of course, $\gamma$ is a random transformation, since it depends on \Xdone{}.)
Then apply the travellers process to $\mgamof{\mXdtwo},\mgamof{\mXdthree}$
to order these:  that is,
        $$(\mYs{2},\mYs{3})=(\mXdtwo,\mXdthree)\quad\mbox{iff }
                \mfthetof{\mgamof{\mXdtwo}} < \mfthetof{\mgamof{\mXdthree}}.$$

Now, having observed \Ys{1}, we know that there are two independent points
below \Ys{1}, arranged according to {\bYthet} (normalized).  Thus the initial 
ranks are
\begin{eqnarray*}
        \lefteqn{\mRs{1}=1,\quad\mRs{2}=2,
                }\\
        \lefteqn{\quad\condprobof{\mRs{3}=2}{\mYs{1},\mYs{2}}
                }\\
        &=&
                1 - \condprobof{\mRs{3}=3}{\mYs{1},\mYs{2}}\\
        &=&1-\mthet
\end{eqnarray*}
so the \sri condition holds.
\end{example}

In the preceding example, {\Xdone} always has the fixed position \Ys{1}, and
the third initial rank \Rs{3} takes the non-extreme value $2$ with positive
probability.  A more
complicated variation is the following:

\begin{example}\label{example:3}
Take $n=5$, and pick two values \inopint{\mthets{1},\mthets{2}}.  Set
        $$\mYs{3}=\mXdthree,$$
and given {\Xdthree}, let
        $$\mIs{2}=\clint{0}{\mXdthree},
                \quad\mIs{1}=\clint{\mXdthree}{1},$$
and set \map{\mgams{i}}{\mIs{i}}{I} to be the affine 
order-preserving bijection for $i=1,2$.

Now, we will \quotes{couple} the other positions as follows
        \begin{equation}\label{eqn:pairs}
                \begin{array}{lcl}
                        \bracepair{\mYs{1}}{\mYs{2}}
                                &=&\bracepair{\mXdowns{4}}{\mXdowns{5}}\\
                        \bracepair{\mYs{4}}{\mYs{5}}
                                &=&\bracepair{\mXdowns{1}}{\mXdowns{2}}
                \end{array}
        \end{equation}
with the specific order within each pair specified by 
$\mfthets{i}\compose\mgams{i}$:  thus,
        $$\pair{\mYs{1}}{\mYs{2}}
                =\pair{\mXdowns{4}}{\mXdowns{5}}
                \quad\mbox{ iff }
                \mfthetsof{1}{\mgamsof{1}{\mXdowns{4}}}
                <
                \mfthetsof{1}{\mgamsof{1}{\mXdowns{5}}}
        $$
(else $\pair{\mYs{1}}{\mYs{2}}=\pair{\mXdowns{5}}{\mXdowns{4}}$) and
        $$\pair{\mYs{4}}{\mYs{5}}
                =\pair{\mXdowns{1}}{\mXdowns{2}}
                \quad\mbox{ iff }
                \mfthetsof{2}{\mgamsof{2}{\mXdowns{1}}}
                <
                \mfthetsof{2}{\mgamsof{2}{\mXdowns{2}}}
        $$
(else $\pair{\mYs{4}}{\mYs{5}}=\pair{\mXdowns{2}}{\mXdowns{1}}$).

Here, by contrast with Example \ref{example:2}, the point {\Xdthree} at the
fixed position {\Ys{3}} is not known until the third observation.  However, 
        $$\condprobof{\mRs{3}=3}{\mYs{1},\mYs{2}}=1$$
and conditioned on any  value of $\mYs{3}=\mXdthree$, we have
        $$
                \condprobof{\mRs{2}=2}{\mYs{1},\mYs{3}}
                =1-
                \condprobof{\mRs{2}=1}{\mYs{1},\mYs{3}}=\mthets{1}
        $$
so that this is also true if we drop the conditioning on \Ys{3}.  Once having observed 
$\mYs{3}=\mXdthree$, we know the next two points lie above \Xdthree, so
        $$\condprobof{\mRs{4}=1}{\mYs{1},\mYs{2},\mYs{3}}=1$$
and independently of $\mYs{1},\mYs{2}$ we know that $\mYs{4},\mYs{5}$ 
satisfy the rank condition
        $$
                \condprobof{\mRs{5}=2}{\mYs{1},\mYs{2},\mYs{3},\mYs{4}}
                =1-
                \condprobof{\mRs{5}=1}{\mYs{1},\mYs{2},\mYs{3},\mYs{4}}
                =\mthets{2}.
        $$
\end{example}

In this example, the rearranged positions were coupled to descending positions
according to a partition into intervals \refeq{pairs}.  However, this
is easily modified:  the reader can check that if for example we couple the
positions via
        \begin{eqnarray*}
        \mYs{2}&=&\mXdthree\\
        \bracepair{\mYs{1}}{\mYs{4}}
                &=&\bracepair{\mXdowns{4}}{\mXdowns{5}}\\
        \bracepair{\mYs{3}}{\mYs{5}}
                &=&\bracepair{\mXdowns{1}}{\mXdowns{2}}
        \end{eqnarray*}
but still use the functions $\mfthets{i}\compose\mgams{i}$ to decide how each 
pair is ordered, then we obtain a rearrangement with
        \begin{eqnarray*}
                \mRs{1}=\mRs{2}&=&\mRs{3}=1\\
                \condprobof{\mRs{4}=3}{\mYs{1},\mYs{2},\mYs{3}}
                &=&1-
                \condprobof{\mRs{4}=4}{\mYs{1},\mYs{2},\mYs{3}}
                =1-\mthets{1}\\
                \condprobof{\mRs{5}=4}{\mYs{1},\mYs{2},\mYs{3},\mYs{4}}
                &=&1-
                \condprobof{\mRs{5}=5}{\mYs{1},\mYs{2},\mYs{3},\mYs{4}}
                =1-\mthets{2}.
        \end{eqnarray*}

One can also increase the number of deterministic positions and/or the
number of positions in any \quotes{coupled} group.

The general construction involves three types of parameters: 
\deffont{fixed positions}, \deffont{switching schemes}, and 
\deffont{jump probabilities}.  Suppose we are working with $n$ variables.

\begin{description}
\item[Fixed positions:] Pick $d$, set $\mns{0}=0$, $\mns{d+1}=n+1$ and
  pick a subsequence            $\mns{1}<\mns{2}<\dots<\mns{d}$
  from \bracelst{1}{n} such that $\mns{i+1}$ is either adjacent
  to $\mns{i}$, or there are at least two intermediate values
  (\ie $\mns{i+1}-\mns{i}\neq 2$).  Let $\mNs{i}\eqdef
  \bracelst{\mns{i-1}+1}{\mns{i}-1}$ (so $\howmany\mNs{i}\geq2$ if
  $\mNs{i}\neq\emptyset$).

\item[Switching schemes:] Pick $d$ distinct positions
  \inblst{\mms{1},\dots,\mms{d}}{1}{n}, and partition           the
  rest of \bracelst{1}{n} into subsets $\mMs{i}$,
  \ilst{i}{1}{d+1}, with $\howmany\mMs{i}=\howmany\mNs{i}$
  for \ilst{i}{1}{d+1}.  Our rearranged positions will be
  coupled to the descending ones via the scheme
  $$\begin{array}{lcll}
    \mYs{\mms{i}}&=&\mXdowns{\mns{i}}\quad&\ilst{i}{1}{d}\\ 
    \setbld{\mYs{j}}{j\in\mMs{i}}                               &=&
    \setbld{\mXdowns{j}}{j\in\mNs{i}}
    \quad&\ilst{i}{1}{d+1}
                \end{array}$$
              \item[Jump probabilities:] For each $i$ such that
                $\mNs{i}\neq\emptyset$, we pick \inopint{\mthets{i}}.
\end{description}

Our map \map{\mu}{\msimplex}{\mSym} defining the rearrangement will then
be defined as follows:  given $\mba=\matuple\in\msimplex$, $\mas{0}=1$,
$\mas{n+1}=0$, $\mmuof{\mba}$ will
satisfy
        \begin{enumerate}
        \item $\mbmus{\mms{i}}=\mns{i}$, \ilst{i}{1}{d};
        \item $j\in\mMs{i}$ iff $\mbmus{j}\in\mNs{i}$;
        \item if $\mNs{i}\neq\emptyset$, let
          $\mIs{i}\eqdef\clint{\mas{\mns{i+1}}}{\mas{\mns{i}}}$,
          take \map{\mgams{i}}{\mIs{i}}{I}
          the affine orientation-preserving bijection,
          and set
          $\mfs{i}\eqdef\mfthets{i}\compose\mgams{i}$;
          then if
          $\mMs{i}=\{\mjs{1}<\mjs{2}<\dots<\mjs{\ell}\}$,
          we define \map{\mu}{\mMs{i}}{\mNs{i}} by the
          condition
          $$\mfsof{i}{\mXdowns{\mmus{\mjs{1}}}}
          <\mfsof{i}{\mXdowns{\mmus{\mjs{2}}}}<\dots
          <\mfsof{i}{\mXdowns{\mmus{\mjs{\ell}}}}.
          $$
        \end{enumerate}

\begin{prop}\label{prop:examples}
        Any rearrangement constructed as above has the \sri property.
\end{prop}

\begin{proof}
  We keep the notation of the construction above.

  First, we determine the initial rank of \Ys{\mms{i}}.  Consider
  $h<\mms{i}$.          Either $h=\mms{p}$ for some $p\neq i$, and
  since         $\mas{\mns{1}}>\mas{\mns{2}}>\dots>\mas{\mns{d}}$,
  $$\mYs{\mms{p}}>\mYs{\mms{i}}\mbox{ iff }p< i$$       or
  $h\in\mMs{p}$ for some $p$, and since
  $\mNs{p}=\setbld{q}{\mns{p-1}<q<\mns{p}}=\setbld{q}{\mas{\mns{p-1}}
    >\mas{q}>\mas{\mns{p}}}$,
  $$\mYs{h}>\mYs{\mms{i}}\mbox{ iff }p\leq i.$$ Thus, by
  \eqref{irank} we have (with probability $1$)
  $\mRs{\mms{i}}=s{}$, with
  \begin{equation}
    \label{eqn:**}
    s=1+\howmany\setbld{j<i}{\mms{j}<\mms{i}}
    +\howmany\left[\bracelst{1}{\mms{i}}\cap \bigcup_{p\leq
        i}\mMs{p}\right].
  \end{equation}

  Now, suppose $k\in\mMs{i}$, and consider $h<k$.  If $h=\mms{p}$ for
  some \inblst{p}{1}{d+1} or  if $h\in\mMs{p}$
  for some $p\neq i$, then we have
  $$\mYs{h}>\mYs{k}\mbox{ iff }p< i.$$  Thus, the only undetermined
  relative sizes are those involving $h\in\mMs{i}$ and
  $h<k$.  Let
  $$r\eqdef\howmany\setbld{h<k}{h\in\mMs{i}}.$$
  Then we know  
  \begin{eqnarray*}
    \lefteqn{\condprobof{\mRs{k}=s}{\indexlst{Y}{1}{k-1},
        \mYs{\mms{i}}, \mYs{\mms{i+1}} } }\\ &&=1-
    \condprobof{\mRs{k}=s+r}{\indexlst{Y}{1}{k-1},
      \mYs{\mms{i}}, \mYs{\mms{i+1}} }\\ 
    &&=1-\mthets{i},
  \end{eqnarray*}
  where $s$ is given by \refeq{**},
  and hence the conditioning on $\mYs{\mms{i}}$ and  $\mYs{\mms{i+1}}$
  can be removed, as in Example \ref{example:3}.
\end{proof}

We pose some unresolved questions concerning the characterization
of rearrangements
with the \sri property.  We use the notation of \refeq{sri} in
Remark \ref{rmk:sri}.

\begin{question}\label{question:char}
        Are rearrangements with the \sri property characterized by the 
        distributions of their rank configurations?
        That is, if {\bY} and $\mbY\pr$ both satisfy \refeq{sri}
        with $\mps{k,\ell}=p\pr_{k,\ell}$ for all $k,\ell$, then does
        it follow that 
                $$\mbY\eqdist\mbY\pr?$$
\end{question}

Theorem \ref{thm:n2} can be viewed as an affirmative answer for $n=2$; a
particular extension would be whether \bYthet{} is characterized by
\begin{eqnarray*}
  \mps{k,1}&=&1-\thet\\
  \mps{k,\ell}&=&0\mbox{ for }1<\ell<k\\
  \mps{k,k}&=&\thet
\end{eqnarray*}
for \ilst{k}{2}{n}.
Two other questions are
raised by our construction above, and by Theorem \ref{thm:indep}.

\begin{question}\label{question:det}
        Is every rearrangement with the \sri property necessarily 
        deterministic?  That is, is it determined by some map 
                $$\map{\mufunct}{\msimplex}{\mSym}?$$
\end{question}

We note that while in our most general examples the ranks do not necessarily 
take extreme values, they still have the property that each rank takes at most
two values.

\begin{question}\label{question:twoval}
        Does there exist a rearrangement with the \sri property for which
        some initial rank can take three or more values with positive
        probability?
\end{question}

The following example shows that strong independence for a single rank, as in
Theorem \ref{thm:bin}, does not alone restrict that rank to two values.

\begin{example}\label{example:4} Take $n=6$, and fix
        $$
                (\mYs{1},\mYs{2},\mYs{3})
                =
                (\mXdowns{1},\mXdowns{3},\mXdowns{5}).
        $$
Then arrange
        $$
                \{\mYs{4},\mYs{5},\mYs{6}\}
                =
                \{\mXdowns{2},\mXdowns{4},\mXdowns{6}\}
        $$
in equiprobable random order.

Now \Rs{4} is equally likely to equal $2, 4$ or $6$, independently of the
values of $\mYs{1},\mYs{2},\mYs{3}$.
\end{example}

Note, of course, that this rearrangement is not deterministic.


\begin{thebibliography}{HK}

\bibitem[GK]{gnedinkrengel} A. V. Gnedin \& U. Krengel, {\em A stochastic game of optimal stopping and order selection}, {\sl Annals Appl.  Prob.} {\bf5}(1995)
310-321.

\bibitem[H1]{hill1}B. Hill, {\em Posterior distribution of percentiles: Bayes' theorem for sampling from  a finite population}, {\sl J. Amer. Stat. Assoc.} {\bf63}(1966), 677-691.

\bibitem[H2]{hill2}$\_\_\_\_$, {\em De Finetti's theorem, induction and $A_{(n)}$ or Bayesian nonparametric predictive inference}, in J. M Bernardo {\em et. al.} (eds.), {\sl Bayesian Statistics 3}. Oxford univ. Press, 1988, 211-241.

\bibitem[HK]{hillkennedy}T. P. Hill \& D. P. Kennedy, {\em Sharp inequalities for optimal stopping with rewards based on ranks}, {\sl Annals Appl. Prob.} {\bf2}(1992), 503-517.


\end{thebibliography}
\end{document}